# Multiple integral representation for functionals of Dirichlet processes

## GIOVANNI PECCATI

*Laboratoire de Statistique Théorique et Appliquée, Université Paris VI, France.*
*E-mail: giovanni.peccati@gmail.com*

We point out that a proper use of the Hoeffding–ANOVA decomposition for symmetric statistics of finite urn sequences, previously introduced by the author, yields a decomposition of the space of square-integrable functionals of a Dirichlet–Ferguson process, written $L^2(D)$, into orthogonal subspaces of multiple integrals of increasing order. This gives an isomorphism between $L^2(D)$ and an appropriate Fock space over a class of deterministic functions. By means of a well-known result due to Blackwell and MacQueen, we show that each element of the $n$th orthogonal space of multiple integrals can be represented as the $L^2$ limit of $U$-statistics with degenerate kernel of degree $n$. General formulae for the decomposition of a given functional are provided in terms of linear combinations of conditioned expectations whose coefficients are explicitly computed. We show that, in simple cases, multiple integrals have a natural representation in terms of Jacobi polynomials. Several connections are established, in particular with Bayesian decision problems, and with some classic formulae concerning the transition densities of multiallele diffusion models, due to Littler and Fackerell, and Griffiths. Our results may also be used to calculate the best approximation of elements of $L^2(D)$ by means of $U$-statistics of finite vectors of exchangeable observations.

## 1. Introduction and preliminaries

Let $(A, \mathcal{A})$ be a Polish space endowed with its Borel $\sigma$-field and consider a finite positive measure $\alpha$ on $(A, \mathcal{A})$. According to [8], given a probability space $(\Omega, \mathcal{F}, \mathbb{P})$, we say that a random probability measure $\{D(C; \omega) : C \in \mathcal{A}\}$, where $\omega \in \Omega$, is a *Dirichlet–Ferguson process* (in the sequel, *DF process*) with parameter $\alpha$ if, for every finite measurable partition $(C_1, \ldots, C_n)$ of $A$, the vector $(D(C_1; \cdot), \ldots, D(C_n; \cdot))$ has a Dirichlet distribution with parameters $(\alpha(C_1), \ldots, \alpha(C_n))$, with the convention that $\alpha(C_i) = 0$ means $D(C_i) = 0$, $\mathbb{P}$-a.s. (throughout the sequel, whenever there is no risk of confusion, we will write $D(C; \omega)$, $D(C; \cdot)$ or $D(C)$ depending on notational convenience). Note that, when $\alpha$ is non-atomic, in the terminology of [24], $D$ is a *normalized gamma process* on $(A, \mathcal{A})$. DF processes were first introduced and analyzed in the fundamental papers [3, 4, 8] and have since played a central role in Bayesian nonparametric statistics (we refer the reader to the above-quoted references, as well as [9, 13] and [20], for basic discussions in this direction; see also [19] for a survey of a large class of random measures related to DF processes).

Now, let $L^2(D) = L^2(D, \mathbb{P})$ denote the Hilbert space of square-integrable functionals of the random measure $D$. The aim of this paper is to obtain an orthogonal decomposition of $L^2(D)$ based on the theory of orthogonal and symmetric $U$-statistics developed in [18] (see also [5]). Such a result is the analogue, for the random measure $D$, of the "chaotic" decompositions of square-integrable functionals of Gaussian processes (see, e.g., [23] and [14] and the references therein) or Lévy processes (see [22] and [15]). In particular, we will show that every element of $L^2(D)$ admits a unique representation as an infinite orthogonal sum of multiple integrals of increasing order with respect to $D$ and therefore that $L^2(D)$ is isomorphic to an appropriate Fock space over a class of deterministic functions. Our results contain as special cases several classic computations contained in [8] and [9], mainly related to Bayesian decision problems. Moreover, they provide an exhaustive characterization of the covariance structure of the elements of $L^2(D)$, for any choice of $(A, \mathcal{A})$ and $\alpha$. In this sense, our results are the infinite-dimensional analogues of the orthogonal polynomial decompositions of functionals of finite Dirichlet vectors, used, for example, by Littler and Fackerell (see [12]) and Griffiths (see [10]) to make explicit the transition density associated with a (finite) multi-allele diffusion model, having the Dirichlet law as stationary measure. Some applications are outlined in Section 1.3 as well as in Sections 6 and 7 below.

To partially illustrate our methods and results in a specific framework, we will first present the example of a simple DF process on $\{0, 1\}$.

## 1.1. Preliminary example: Beta random variables and Jacobi polynomials

Fix real numbers $\alpha_1, \alpha_0 > 0$ and consider a Beta random variable $\eta(\omega)$ with values in $[0, 1]$ and parameters $(\alpha_1, \alpha_0)$. This means that, for every Borel set $C$,

$$\mathbb{P}(\eta \in C) = \frac{1}{B(\alpha_1, \alpha_0)} \int_{C \cap [0,1]} x^{\alpha_1 - 1}(1-x)^{\alpha_0 - 1} \, \mathrm{d}x, \tag{1}$$

where $B(\cdot, \cdot)$ is the Beta function, defined as $B(s, t) = \int_0^1 x^{s-1}(1-x)^{t-1} \, \mathrm{d}x$ (see, e.g., [1]). We may interpret $\eta$ as a random parameter, determining a random probability measure $D(\cdot, \omega)$ on $\{0, 1\}$, via the relations

$$D(\{1\}, \omega) = \eta(\omega) = 1 - D(\{0\}, \omega). \tag{2}$$

The measure $D(\cdot, \omega)$, as defined in (2), is the most elementary example of a DF process and corresponds, in particular, to the case $A = \{0, 1\}$ and $\alpha(\cdot) = \alpha_1 \delta_1(\cdot) + \alpha_0 \delta_0(\cdot)$, where $\delta_x$ stands for the Dirac measure concentrated at $x$. To simplify, we adopt the notation $p_{\alpha_1,\alpha_0}(x) = B(\alpha_1,\alpha_0)^{-1}x^{\alpha_1-1}(1-x)^{\alpha_0-1}$ and define $L^2(\eta)$ and $\mathcal{P}_n(\eta)$, $n \geq 0$, to be, respectively, the space of square-integrable functionals of $\eta$ and the subspace of $L^2(\eta)$ composed of random variables of the form $\pi_n(\eta)$, where $\pi_n(\cdot)$ is a polynomial of order $n$. Note that $L^2(\eta) = L^2(D)$, $\mathcal{P}_0(\eta) = \Re$, $\mathcal{P}_n(\eta) \subset \mathcal{P}_{n+1}(\eta)$ and the union of the $\mathcal{P}_n(\eta)$'s is total in $L^2(\eta)$; we also set $\mathcal{J}_0(\eta) := \Re$ and $\mathcal{J}_n(\eta) := \mathcal{P}_n(\eta) \cap \mathcal{P}_{n-1}(\eta)^{\perp}$, $n \geq 1$, where "$\perp$" stands for the orthogonality relation in $L^2(\eta)$. It is well known that the orthogonal sequence of subspaces $\{\mathcal{J}_n(\eta) : n \geq 0\}$ can be exhaustively characterized in terms of *Jacobi polynomials* (again, see [1], Section 22), defined, for $n \geq 0$, $q > 0$ and $p > q - 1$, as $G_n(p,q,x) := \sum_{a=0,\ldots,n} g_{n,a}(p,q)x^a$, where $g_{n,a}(p,q) := \binom{n}{a}(-1)^{n-a}\frac{\Gamma(q+n)\Gamma(p+a+n)}{\Gamma(p+2n)\Gamma(a+q)}$. Indeed, for $\alpha_1$ and $\alpha_0$ as before, one can prove that the sequence of *modified Jacobi polynomials*, defined through the relation

$$J_n^{\alpha_1,\alpha_0}(x) = G_n(\alpha_1 + \alpha_0 - 1, \alpha_1, x)\sqrt{k_n(\alpha_1,\alpha_0)}$$
$$= \sum_{a=0}^{n} c_{n,a}(\alpha_1 + \alpha_0 - 1, \alpha_1)x^a, \qquad n \geq 0, \tag{3}$$

where

$$k_n(\alpha_1,\alpha_0) = \frac{(2n+\alpha_1+\alpha_0-1)\Gamma^2(2n+\alpha_0+\alpha_1-1)}{n!\Gamma(n+\alpha_1)\Gamma(n+\alpha_0)\Gamma(n+\alpha_1+\alpha_0-1)}B(\alpha_1,\alpha_0),$$

$$c_{n,a}(\alpha_1+\alpha_0-1,\alpha_1) = \sqrt{k_n(\alpha_1,\alpha_0)}g_{n,a}(\alpha_1+\alpha_0-1,\alpha_1),$$

is such that $\int_0^1 J_n^{\alpha_1,\alpha_0}(x)J_m^{\alpha_1,\alpha_0}(x)p_{\alpha_1,\alpha_0}(x)\,\mathrm{d}x = 0$ or $1$, according to whether $m \neq n$ or $m = n$, thus implying that the class $\{J_n^{\alpha_1,\alpha_0} : n \geq 0\}$ is a family of *orthogonal polynomials* associated with the weight function $p_{\alpha_1,\alpha_0}$ on the interval $[0,1]$. This immediately yields that, for every $n \geq 0$, $X \in \mathcal{J}_n(\eta)$ if and only if $X = cJ_n^{\alpha_1,\alpha_0}(\eta)$ for some real constant $c$ and therefore that every $F \in L^2(\eta)$ admits a unique representation of the form

$$F = \mathbb{E}(F) + \sum_{n=1}^{\infty} c_n J_n^{\alpha_1,\alpha_0}(\eta), \tag{4}$$

where the real constants $c_n$ are such that $\sum c_n^2 < +\infty$ (i.e., the series on the right-hand side of (4) converges in $L^2(\eta)$). It is not difficult to see (see Section 5 below for a complete discussion of this point) that, for every $n$, the random variable $c_n J_n^{\alpha_1,\alpha_0}(\eta)$ can be (uniquely) written in the form $\int_{\{0,1\}^n} \phi_n \, \mathrm{d}D^{\otimes n}$, where $\mathrm{d}D^{\otimes n}$ is the random product measure on $\{0,1\}^n$ generated by the random probability defined in (2) and $\phi_n$ is a well-chosen *symmetric* kernel on $\{0,1\}^n$. This implies, in particular, that every $F \in L^2(\eta)$ admits a decomposition as an infinite orthogonal sum of multiple random integrals, that is,

$$F = \mathbb{E}(F) + \sum_{n=1}^{\infty} \int_{\{0,1\}^n} \phi_n \, \mathrm{d}D^{\otimes n}. \tag{5}$$

We will complete the example above in Section 5 by showing that the kernels $\phi_n$ have a natural interpretation in terms of $U$-statistics. The connections between our results and other special polynomials in several variables are discussed in Section 6. Before that, we shall generalize the representations (4) and (5) by obtaining an analogous orthogonal decomposition of the space $L^2(D)$ associated with a DF process $D$, with an arbitrary parameter $\alpha(\cdot)$ and defined on a general Polish space $(A, \mathcal{A})$.

### 1.2. Discussion of the main results

Let $\alpha$ be a finite measure on the Polish space $(A, \mathcal{A})$ and let $D$ be a DF process of parameter $\alpha$. To obtain our main results – and to be able to use the theory developed in [18] – we shall suppose that the law of $D$ is the *de Finetti measure*, that is, that $D$ is the *directing measure* of an infinite exchangeable sequence $\mathbf{X} = \{X_n : n \geq 1\}$ of random variables with values in $(A, \mathcal{A})$. This means that the sequence $\mathbf{X}$ is defined on the same probability space as $D$ and that, conditioned on $D$, $\mathbf{X}$ is composed of i.i.d. random variables with common law equal to $D$ (see [2] for an exhaustive discussion of this point). Note that, given a general random probability measure $M(\cdot; \omega)$, there always exists (on a possibly enlarged probability space) an exchangeable sequence $\mathbf{Y}$ such that $M$ is the directing measure of $\mathbf{Y}$. Then, according to, for example, [4], $\mathbf{X}$ must be an infinite *generalized Pólya urn sequence* with parameter $\alpha$, as defined in the above-quoted reference and in Section 2 below. Note that, in this case, $D$ is automatically the a.s. limit of the sequence of empirical measures generated by $\mathbf{X}$.

The principal achievement of the present paper is to prove (Theorem 1) that every $F \in L^2(D)$ has a unique representation of the type

$$
\begin{aligned}
F &= \mathbb{E}(F) + \sum_{n \geq 1} \int_{A^n} h_{(F,n)}(a_1, \ldots, a_n) D^{\otimes n}(\mathrm{d}a_1, \ldots, \mathrm{d}a_n) \\
  &= \mathbb{E}(F) + \sum_{n \geq 1} \int_{A^n} h_{(F,n)} \, \mathrm{d}D^{\otimes n},
\end{aligned}
\tag{6}
$$

where $D^{\otimes n}$ indicates the $n$-dimensional (random) product measure associated with $D$, the series converges in $L^2$ and the kernels $h_{(F,n)}$, $n \geq 1$, are deterministic, symmetric and such that, for every $n$,

$$
\mathbb{E}(h_{(F,n)}(\mathbf{X}_n)^2) < +\infty \quad \text{and} \quad \mathbb{E}(h_{(F,n)}(\mathbf{X}_n) \mid \mathbf{X}_{n-1}) = 0, \qquad \mathbb{P}\text{-a.s.} \tag{7}
$$

Here, $\mathbf{X}_n = (X_1, \ldots, X_n)$ represents, for every $n \geq 1$, the first $n$ instants of the Pólya sequence $\mathbf{X}$ introduced at the beginning of this subsection. Consistent with the notation of [18] and [16], Chapters 9 and 10, and for $n \geq 1$, the class of symmetric functions $h$ on $A^n$ satisfying condition (7) is denoted $\Xi_n(\mathbf{X})$. The functions $h_{(F,n)} \in \Xi_n(\mathbf{X})$ appearing in (6) may be interpreted, for every $n$, as completely degenerate kernels of symmetric $U$-statistics (see, e.g., [11]) based on a truncation of the sequence $\mathbf{X}$. As a consequence (see Proposition 3 below), an application of the results contained in [18] yields that the sequence $\int_{A^n} h_{(F,n)} \, \mathrm{d}D^{\otimes n}$, $n \geq 1$, appearing in (6) enjoys the following isometric property: for every $n, m \geq 1$,

$$\mathbb{E}\left(\int_{A^n} h_{(F,n)} \, \mathrm{d}D^{\otimes n} \int_{A^m} h_{(F,m)} \, \mathrm{d}D^{\otimes m}\right) = \epsilon_{m,n} \times c(n, \alpha(A)) \mathbb{E}(h_{(F,n)}(\mathbf{X}_n)^2), \quad (8)$$

where $c(n, \alpha(A)) := \prod_{l=1}^{n}(n - l + 1)/(\alpha(A) + n + l - 1)$ and $\epsilon_{m,n}$ equals 0 or 1 according to whether $m \neq n$ or $m = n$. As anticipated (see Proposition 5 below), a random variable of the type $\int_{A^n} h \, \mathrm{d}D^{\otimes n}$, $h \in \Xi_n(\mathbf{X})$, represents the infinite-dimensional analogue of the modified Jacobi polynomial introduced in (3). Note that formula (8) determines an isomorphism between $L^2(D)$ and the orthogonal sum

$$\bigoplus_{n \geq 0} \sqrt{c(n, \alpha(A))} \Xi_n(\mathbf{X}) \simeq \bigoplus_{n \geq 0} \sqrt{c(n, \alpha(A))} SH_n(\mathbf{X}_n), \quad (9)$$

where "$\simeq$" indicates a Hilbert space isomorphism, $SH_0 = \Re$ and $SH_n(\mathbf{X}_n)$ is the $n$th *symmetric Hoeffding space* associated with the finite Pólya urn sequence $\mathbf{X}_n$ (see [18], Section 3, and [17], as well as Section 3 below). More to the point, a recursive formula is given (Theorem 2) to explicitly calculate real coefficients $\{\theta^{(n,k)} : n \geq 1, 1 \leq k \leq n\}$ that depend uniquely on $\alpha(A)$ and satisfy the relation

$$h_{(F,n)}(a_1, \ldots, a_n) = \sum_{k=1}^{n} \theta^{(n,k)} \sum_{1 \leq j_1 < \cdots < j_k \leq n} \mathbb{E}(F - \mathbb{E}(F) \mid X_1 = a_{j_1}, \ldots, X_k = a_{j_k}) \quad (10)$$

for every $F \in L^2(D)$. It is worth noting that (10) is quite explicit since, according to, for example, [8], Theorem 1, for every $k \geq 1$ and every $(a_1, \ldots, a_k) \in A^k$, the law of $D$ under the conditioned measure $\mathbb{P}(\cdot \mid X_1 = a_1, \ldots, X_k = a_k)$ is that of a DF process with parameter $\alpha + \sum_{i=1}^{k} \delta_{a_i}$, where $\delta_a$ indicates the Dirac mass at $a$. Such a stability property of the class of DF processes is usually summarized by saying that DF processes are *conjugate* (see [20]). Also, observe that, unlike other chaotic decompositions, to obtain the explicit formula (10), we do not need any regularity assumption on $F$ (see, e.g., [23] for Wiener chaos, where the regularity assumptions are related to weak differentiability, in the sense of Shigekawa–Malliavin).

## 1.3. Some motivations from Bayesian nonparametric statistics (Bayesian estimation of conditional variances)

Our results contain, as special cases, several computations from [8], Sections 4 and 5, and therefore have some immediate applications to (nonparametric) Bayesian statistical decision problems. As an illustration, consider the following setup (see [8] for further details). We are given a sequence of random variables $\mathbf{X} = \{X_i : i \geq 0\}$, modeling the observations of a random phenomenon with values in $(A, \mathcal{A})$, such that $X_0 = x_0 \in A$ and, conditionally on the realization of a random probability measure $D$, the sequence $\{X_i : i \geq$ 1} is i.i.d. with common law equal to $D$. We shall use the notation $(X_0, X_1, \ldots, X_n) = \mathbf{X}_n$, $n \geq 0$, and suppose that $D$ has the law (known as the *prior distribution*) of a DF process with parameter $\alpha$, with $\alpha(A) < +\infty$. In particular, the measure $\alpha$ (which completely determines the law of $D$) is a mathematical representation of the *initial information* of the observer. We also consider a functional $F \in L^2(D)$ with the form of a *conditional variance*, that is,

$$F = V(h) = \int_{A^\infty} h(\mathbf{a})^2 D^\infty(\mathrm{d}\mathbf{a}) - \left( \int_{A^\infty} h(\mathbf{a}) D^\infty(\mathrm{d}\mathbf{a}) \right)^2$$

$$= \mathbb{E}[(h(\mathbf{X}) - \mathbb{E}[h(\mathbf{X}) \mid D])^2 \mid D],$$

where $h : A^\infty \mapsto \Re$ is such that $\mathbb{E}(h(\mathbf{X})^2) < +\infty$ and $D^\infty$ is the canonical (infinite) product measure induced by $D$ on $A^\infty$. Note that, except in trivial cases, $F$ is a function of the whole (infinite) sequence $\mathbf{X}$: it follows that its value cannot be inferred from *any* finite set of observations. Given the sample $\mathbf{x}_n = (x_0, x_1, \ldots, x_n)$ of the first $n+1$ observations ($n \geq 0$), we therefore face the following decision problem: provide an estimation of $F$ by choosing the square-integrable statistic $\widehat{h}_{V(h)}(\mathbf{x}_n)$ that minimizes the (conditional) expected square loss

$$L(\widehat{h}_{V(h)}; \alpha; \mathbf{x}_n) = \mathbb{E}[(F - \widehat{h}_{V(h)}(\mathbf{X}_n))^2 \mid (X_0, \ldots, X_n) = \mathbf{x}_n]; \tag{11}$$

observe that, with this notation, for $n = 0$,

$$L(\widehat{h}_{V(h)}; \alpha; \mathbf{x}_0) = L(\widehat{h}_{V(h)}; \alpha; x_0) = \mathbb{E}[(V(h) - \widehat{h}_{V(h)}(x_0))^2].$$

Since the conjugacy of DF processes implies that, under the probability $\mathbb{P}[\cdot \mid (X_0, \ldots, X_n) = \mathbf{x}_n]$, $D$ has the law of a DF process with parameter $\alpha_{\mathbf{x}_n} := \alpha + \sum_{j=1}^{n} \delta_{x_j}$, for any choice of $\alpha$ and $\mathbf{x}_n$, we have $L(\widehat{h}_{V(h)}; \alpha; \mathbf{x}_n) = L(\widehat{h}_{V(h)}; \alpha_{\mathbf{x}_n}; x_0)$, with $\delta_x$ the Dirac mass in $x$. Elementary computations therefore yield

$$\widehat{h}_{V(h)}(\mathbf{x}_n) = \mathbb{E}\left[ \int_{A^\infty} h(\mathbf{a})^2 \widetilde{D}^\infty(\mathrm{d}\mathbf{a}) - \left( \int_{A^\infty} h(\mathbf{a}) \widetilde{D}^\infty(\mathrm{d}\mathbf{a}) \right)^2 \right]$$

$$= \mathbb{E}[h(\mathbf{X})^2 \mid (X_0, \ldots, X_n) = \mathbf{x}_n] - \mathbb{E}[V(h)^2 \mid (X_0, \ldots, X_n) = \mathbf{x}_n], \tag{12}$$

where $\widetilde{D}$ is a DF process with parameter $\alpha_{\mathbf{x}_n}$. Now, consider the following decomposition of $V(h)$ under the measure $\mathbb{P}[\cdot \mid (X_0, \ldots, X_n) = \mathbf{x}_n]$:

$$V(h) = \mathbb{E}[h(\mathbf{X}) \mid (X_0, \ldots, X_n) = \mathbf{x}_n] + \sum_{k=1}^{\infty} \int_{A^k} h_{(V(h),k)}^{(\mathbf{x}_n)}(a_1, \ldots, a_k) \, \mathrm{d}\widetilde{D}^{\otimes k}(a_1, \ldots, a_k),$$

where the kernels $h_{(V(h),k)}^{(\mathbf{x}_n)}$ can be obtained by applying formula (10) to a DF process with parameter $\alpha_{\mathbf{x}_n}$. We can conclude from (12) and the orthogonality of the $h_{(V(h),k)}^{(\mathbf{x}_n)}$

that

$$
\widehat{h}_{V(h)}(\mathbf{x}_n) = \mathbf{Var}[h(\mathbf{X}) \mid (X_0, \ldots, X_n) = \mathbf{x}_n]
$$

$$
- \sum_{k=1}^{\infty} c(k, \alpha(A) + n) \times \mathbb{E}[h_{(V(h),k)}^{(\mathbf{x}_n)}(X_{n+1}, \ldots, X_{n+k})^2 \mid (X_0, \ldots, X_n) = \mathbf{x}_n], \tag{13}
$$

where $\mathbf{Var}[\cdot \mid (X_0, \ldots, X_n) = \mathbf{x}_n]$ stands for the conditional variance. We stress that the kernels $h_{(V(h),k)}^{(\mathbf{x}_n)}$ in (13) are explicitly known, due to formula (10). Moreover, it will become clear from the subsequent analysis that if $V(h) = \mathbb{E}[(h(\mathbf{X}_m) - \mathbb{E}[h(\mathbf{X}_m) \mid D])^2 \mid D]$ for some $1 \leq m < +\infty$, then $h_{(V(h),k)}^{(\mathbf{x}_n)} = 0$ for $k > m$ and therefore the right-hand side of (13) is just a finite sum. In particular, formula (13) generalizes the computations contained in [8], Section 5(e). For instance, if $h(\mathbf{a}) = h(a_1)$ (so that $V(h) = \int_A h^2 \, \mathrm{d}D - (\int_A h \, \mathrm{d}D)^2$), then (13) reduces to the well-known formula (see, e.g., [8], page 226)

$$
\widehat{h}_{V(h)}(\mathbf{x}_n) = \frac{\alpha(A) + n}{\alpha(A) + n + 1} \mathbf{Var}[h(X_{n+1}) \mid (X_0, \ldots, X_n) = \mathbf{x}_n].
$$

## 1.4. Further remarks and organization of the paper

As discussed below, Theorems 1 and 2 represent a logical continuation of the results contained in [18], Section 5, where we obtained the explicit Hoeffding–ANOVA decomposition for symmetric statistics of vectors of exchangeable observations that are (finite) Generalized Urn Sequences (GUS). The class of GUS contains, as special cases, vectors of i.i.d. random variables, as well as extractions without replacement from a finite population and truncated Pólya urn sequences. The results of this paper are mainly obtained by properly extending to the infinite-dimensional case the content of the above-quoted reference.

The analysis contained in this work is also related to another statistical problem. Supposing that we are given a vector $(X_1, \ldots, X_n)$ of exchangeable observations that are the first $n$ instants of the infinite Pólya sequence $\mathbf{X}$, which is the best approximation of a generic element of $L^2(D)$ by means of $U$-statistics that are based exclusively on $(X_1, \ldots, X_n)$ and how can one compute the corresponding quadratic error? Of course, one can ask the same question for a general infinite exchangeable sequence and for any square-integrable functional of its directing measure. In the last section of the paper, we shall show that a general solution to this problem is contained in formulae (8) and (10) above, as well as in the calculations performed in [18].

The paper is organized as follows. In the next section, we recall some results about Dirichlet processes and exchangeable sequences which are due to Blackwell and Mac-Queen. In Section 3, some preliminaries about urn sequences and Hoeffding–ANOVA decompositions are presented. Section 4 contains the statements and proofs of the two main theorems of this work. In Section 5, we complete the study of the simple Dirichlet process introduced in (2) and establish, in this case, an explicit relation between Jacobi polynomials and the multiple random integrals appearing in (6). In Section 6, several connections are discussed between our decomposition of the space $L^2(D)$ and the family of generalized Appell–Jacobi polynomials used in [12] and [10] to make explicit the transition density of a Wright–Fisher diffusion process. Section 7 contains further applications and examples.

The results of this paper have been partially announced in [17].

## 2. Blackwell–MacQueen construction of the Dirichlet process

The main idea of [4] is that a general Dirichlet process can be represented as the limit of the empirical measures associated with an infinite, exchangeable sequence of observations and that such a sequence can be taken to be a generalization of the so-called *Pólya urn scheme*. The notation of the previous section is maintained throughout the sequel.

Suppose that a sequence of random variables $\mathbf{X} = \{X_n : n \geq 1\}$ is defined on the probability space $(\Omega, \mathcal{F}, \mathbb{P})$, taking values in $(A, \mathcal{A})$ and such that, for every $k \geq 1$ and every $1 \leq j_1 < \cdots < j_k < +\infty$,

$$\mathbb{P}(X_{j_1} \in da_1, \ldots, X_{j_k} \in da_k) = \prod_{i=1}^{k} \frac{\alpha(da_i) + \sum_{l=1}^{i-1} \delta_{a_l}(da_i)}{\alpha(A) + i - 1}. \tag{14}$$

We can think of $\mathbf{X}$ as an infinite sequence of extractions from an urn $A$, whose initial composition is given by the measure $\alpha$, according to the classic *Pólya scheme*: at each step, a ball is extracted and two balls of the same color are placed in the urn before the next extraction. It is also clear, due to (14), that $\mathbf{X}$ is an infinite exchangeable sequence, in the sense that its law is invariant under finite permutations of the index set $\{1, 2, \ldots\}$. We call $\mathbf{X}$ an (infinite) *Pólya sequence* with parameter $\alpha$. In the terminology of [18], Section 5, we have that, for each $N \geq 2$, the vector $\mathbf{X}_N = (X_1, \ldots, X_N)$ is a Generalized Urn Sequence (GUS) of length $N$ with parameters $\alpha$ and $c = 1$. This implies that, for every $k \geq 1$ and every $(a_1, \ldots, a_k) \in A^k$, under the conditioned probability $\mathbb{P}(\cdot \mid X_1 = a_1, \ldots, X_k = a_k)$, the sequence $\{X_{k+n} : n \geq 1\}$ is an infinite Pólya sequence with parameter $\alpha + \sum_{l=1}^{k} \delta_{a_l}$. The next result, which is proved in [4], states that the sequence of empirical measures generated by $\mathbf{X}$ converges almost surely to a DF process with parameter $\alpha$ and that the law of such a process coincides with the de Finetti measure associated with $\mathbf{X}$.

**Theorem (Blackwell and MacQueen [4]).** *Using the previous notation for every $n$, define*

$$\mathbb{P}_n(C; \omega) = \frac{1}{n} \sum_{i=1}^{n} \mathbf{1}_C(X_i(\omega)), \qquad C \in \mathcal{A},$$

*to be the empirical measure associated with the vector $\mathbf{X}_n$. Then,*

(a) *as $n$ goes to infinity, the random measure $\mathbb{P}_n(\cdot; \omega)$ converges $\mathbb{P}$-a.s. to a random discrete probability $D(\cdot; \omega)$ on $(A, \mathcal{A})$;*

(b) *the measure $D$ appearing in* (a) *is a DF process with parameter $\alpha$;*

(c) *given $D$, the variables $X_1, X_2, \ldots$ composing the sequence $\mathbf{X}$ are independent and identically distributed with law $D$.*

By inspection of the proof contained in [4], the convergence in item (a) of the previous theorem can be interpreted in the following sense: there exists $\Omega_* \in \mathcal{F}$ such that $\mathbb{P}(\Omega_*) = 1$ and, for every $\omega \in \Omega_*$,

$$\mathbb{P}_n(C; \omega) \longrightarrow D(C; \omega) \qquad \forall C \in \mathcal{A}.$$

This implies that, almost surely, $\mathbb{P}_n$ weakly converges to $D$. The reader is also referred to [19], where it is shown that weak convergence may be replaced by convergence in total variation. Also, note that dominated convergence implies that, for any measurable set $C$, $\mathbb{E}[(\mathbb{P}_n(C) - D(C))^2] \to 0$. In the classical terminology of [8], item (c) of the previous theorem states that, for every $k \geq 1$ and every $j_1 \neq \cdots \neq j_k$, the vector $(X_{j_1}, \ldots, X_{j_k})$ is a *sample of size $k$* from $D$, where $D$ is a DF process appearing as the limit of the sequence $\mathbb{P}_n$, $n \geq 1$. From now on, when considering a DF process $D$ with parameter $\alpha$, we will always assume that such a process is the a.s. limit of the sequence $\mathbb{P}_n$ associated with a Pólya sequence $\mathbf{X}$ with the same parameter.

# 3. Urn sequences and Hoeffding–ANOVA decompositions

In this section, we recall the results of [18] and [16] that are related to generalized Pólya urn sequences. We start by introducing some notation, mostly borrowed from the above-quoted references.

Fix $N \geq 1$. For any $n \in \{0, 1, \ldots, N\}$, we define

$$V_N(n) := \{\mathbf{k}_{(n)} = (k_1, \ldots, k_n) : 1 \leq k_1 < \cdots < k_n \leq N\},$$

where $\mathbf{k}_{(0)} := 0$ and $V_N(0) = \{0\}$, and also $V_\infty(n) = \bigcup_{N \geq 1} V_N(n)$. For $n \geq m \geq 1$, $\mathbf{l}_{(m)} \in V_\infty(m)$ and $\mathbf{k}_{(n)} \in V_\infty(n)$, $\mathbf{l}_{(m)} \wedge \mathbf{k}_{(n)}$ is the set $\{l_i : l_i = k_j$ for some $j = 1, \ldots, n\}$ written as an element of $V_\infty(r)$, where $r := \mathrm{Card}\{\mathbf{l}_{(m)} \wedge \mathbf{k}_{(n)}\}$. Analogously, for any $n, m \geq 0$, $\mathbf{k}_{(n)} \backslash \mathbf{l}_{(m)}$ denotes the set $\mathbf{k}_{(n)} \cap (\mathbf{l}_{(m)})^c$ written as an element of the class $V_\infty(n - r)$. Finally, given $\mathbf{k}_{(n)} \in V_\infty(n)$ and a vector $\mathbf{h}_{(m)} = (h_1, \ldots, h_m)$, by $\mathbf{h}_{(m)} \subset \mathbf{k}_{(n)}$, we mean that $\mathbf{h}_{(m)} \in V_\infty(m)$ and that for every $i \in \{1, \ldots, m\}$, there exists $j \in \{1, \ldots, n\}$ such that $k_j = h_i$.

Now, consider a Pólya sequence $\mathbf{X}$ such as the one defined in the previous section. For any $n \geq 0$ and every $\mathbf{j}_{(n)} \in V_\infty(n)$, we write $\mathbf{X}_{\mathbf{j}_{(n)}} = (X_{j_1}, \ldots, X_{j_n})$, with $\mathbf{X}_0 = 0$. We recall that the exchangeability of $\mathbf{X}$ implies that, for any $0 \leq r \leq m \leq n$ and for any symmetric statistic $T$ on $A^n$ such that $\mathbb{E}[|T(\mathbf{X}_n)|] < +\infty$, there exists a function $[T]_{n,m}^{(r)}$ on $A^m$, symmetric in the first $r$ variables and in the last $m - r$ variables such that, for every $\mathbf{j}_{(n)} \in V_\infty(n)$ and $\mathbf{i}_{(m)} \in V_\infty(m)$ satisfying $\mathrm{Card}(\mathbf{i}_{(m)} \wedge \mathbf{j}_{(n)}) = r$,

$$\mathbb{E}[T(\mathbf{X}_{\mathbf{j}_{(n)}}) \mid \mathbf{X}_{\mathbf{i}_{(m)}}] = [T]^{(r)}_{n,m}(\mathbf{X}_{\mathbf{j}_{(n)} \wedge \mathbf{i}_{(m)}}, \mathbf{X}_{\mathbf{j}_{(n)} \setminus \mathbf{i}_{(m)}}), \qquad \text{a.s.-}\mathbb{P}.$$

In [18], we have provided a complete characterization of the *symmetric Hoeffding spaces* associated with the random vector $\mathbf{X}_{\mathbf{j}_{(N)}}$ for any $N \geq 1$ and any $\mathbf{j}_{(N)} \in V_\infty(N)$. More precisely, we start by writing $L^2(\mathbf{X}_{\mathbf{j}_{(N)}})$ for the Hilbert space of real-valued and square-integrable functionals of $\mathbf{X}_{\mathbf{j}_{(N)}}$ and $L^2_s(\mathbf{X}_{\mathbf{j}_{(N)}})$ for the subspace of $L^2(\mathbf{X}_{\mathbf{j}_{(N)}})$ composed of symmetric functionals. Of course, for every $\mathbf{j}_{(N)}$, the space $L^2_s(\mathbf{X}_{\mathbf{j}_{(N)}})$ coincides with the space of square-integrable functionals of the empirical random measure generated by the vector $\mathbf{X}_{\mathbf{j}_{(N)}}$. We eventually set $L^2(\mathbf{X})$ to be the space of square-integrable functionals of the sequence $\mathbf{X}$ (note that $L^2(D) \subset L^2(\mathbf{X})$).

According to [18], the collection $\{SH_i(\mathbf{X}_{\mathbf{j}_{(N)}}), i = 0, \ldots, N\}$ of *symmetric Hoeffding spaces* associated with $\mathbf{X}_{\mathbf{j}_{(N)}}$ is defined as follows. Let $SU_0(\mathbf{X}_{\mathbf{j}_{(N)}}) := \Re$ and, for $i = 1, \ldots, N$,

$$SU_i(\mathbf{X}_{\mathbf{j}_{(N)}}) := \overline{\text{v.s.} \left\{ T : T = \sum_{\mathbf{j}_{(i)} \subset \mathbf{j}_{(N)}} g(\mathbf{X}_{\mathbf{j}_{(i)}}), g(\mathbf{X}_{\mathbf{j}_{(i)}}) \in L^2_s(\mathbf{X}_{\mathbf{j}_{(i)}}) \right\}}^{L^2(\mathbf{X})},$$

where v.s. $\{C\}$ is the minimal vector space containing $C$ and

$$SH_0(\mathbf{X}_{\mathbf{j}_{(N)}}), = SU_0(\mathbf{X}_{\mathbf{j}_{(N)}}),$$
$$SH_i(\mathbf{X}_{\mathbf{j}_{(N)}}) = SU_i(\mathbf{X}_{\mathbf{j}_{(N)}}) \ominus SU_{i-1}(\mathbf{X}_{\mathbf{j}_{(N)}}), \qquad i = 1, \ldots, N,$$

where $\ominus$ denotes orthogonal difference between Hilbert spaces. The reader is referred to [18] and the references therein for more details about the use and interpretation of Hoeffding spaces. As discussed in [18], $L^2_s(\mathbf{X}_{\mathbf{j}_{(N)}})$ differs from $SU_i(\mathbf{X}_{\mathbf{j}_{(N)}})$ for every $i = 0, \ldots, N-1$. This means, in particular, that

$$SU_i(\mathbf{X}_{\mathbf{j}_{(N)}}) = \bigoplus_{a \leq i} SH_a(\mathbf{X}_{\mathbf{j}_{(N)}}) \subsetneq L^2_s(\mathbf{X}_{\mathbf{j}_{(N)}}) = SU_N(\mathbf{X}_{\mathbf{j}_{(N)}}) = \bigoplus_{a=0,\ldots,n} SH_a(\mathbf{X}_{\mathbf{j}_{(N)}}).$$

Now, consider the measure $\alpha$ on $(A, \mathcal{A})$ that determines the law of $\mathbf{X}$. We shall use the following real constants:

$$\Phi(n, m, r, p) := (m - r)_{(m-r-p)} \frac{\prod_{s=1}^{m-(r+p)}[\alpha(A) + r + p + s - 1]}{\prod_{s=1}^{m-r}[\alpha(A) + n + s - 1]}, \tag{15}$$

where $1 \leq m \leq n$, $0 \leq r \leq m$, $0 \leq p \leq m - r$, $\alpha(A) + n + m - r > 0$, $(a)_{(b)} := a!/b!$ for $a \geq b$ and $\prod_{s=1}^0 = 1 = 0^0$, by definition, and, for $1 \leq q \leq m \leq n \leq N$,

$$\Psi_N(q, n, m) := \sum_{r=0}^q \binom{q}{r} \binom{N-n}{m-r}_* \Phi(n, m, r, q - r) \tag{16}$$

with $\binom{a}{b}_* := \binom{a}{b}\mathbf{1}_{(a \geq b)}$. We also introduce the coefficients

$$\{\theta_N^{(k,a)} : 1 \leq k \leq N, 1 \leq a \leq k\}$$

that are recursively defined by the set of conditions $\{\mathbf{S}_N(k), k = 1, \ldots, N-1\}$ given by

$$\mathbf{S}_N(k) := \begin{cases} \theta_N^{(k,k)} = \Psi_N(k,k,k)^{-1}, \\ \sum_{i=q}^{k} \sum_{j=q}^{i} \theta_N^{(i,j)} \Psi_N(q,k,j) = 0, \qquad q = 1, \ldots, k-1, \end{cases} \tag{17}$$

and $\theta_N^{(N,a)} := -\sum_{s=a}^{N-1} \theta_N^{(s,a)}$ for $a = 1, \ldots, N-1$ and $\theta_N^{(N,N)} = \Psi_N(N,N,N)^{-1} = 1$. We further set

$$\theta_{N*}^{(k,a)} := \theta_N^{(k,a)} \binom{N-a}{k-a}^{-1}.$$

The following proposition stems from the main results obtained [18], Section 5. It provides an algorithm to project symmetric statistics onto Hoeffding spaces.

**Proposition 1 (Peccati [18]).** *Using the above notation and assumptions, fix* $\mathbf{j}_{(N)} \in V_\infty(N)$ *and let* $T$ *be a centered element of* $L_s^2(\mathbf{X}_{\mathbf{j}_{(N)}})$. *Then, for* $s = 1, \ldots, N$,

$$\pi[T, SH_s](\mathbf{X}_{\mathbf{j}_{(N)}}) = \sum_{\mathbf{j}_{(s)} \subset \mathbf{j}_{(N)}} \left[ \sum_{a=1}^s \theta_{N*}^{(s,a)} \sum_{\mathbf{j}_{(a)} \subset \mathbf{j}_{(s)}} [T]_{N,a}^{(a)}(\mathbf{X}_{\mathbf{j}_{(a)}}) \right] = \sum_{\mathbf{j}_{(s)} \subset \mathbf{j}_{(N)}} \phi_T^{(s)}(\mathbf{X}_{\mathbf{j}_{(s)}}),$$

*where*

$$\phi_T^{(s)}(\mathbf{X}_{\mathbf{j}_{(s)}}) = \left[ \sum_{a=1}^s \theta_{N*}^{(s,a)} \sum_{\mathbf{j}_{(a)} \subset \mathbf{j}_{(s)}} [T]_{N,a}^{(a)}(\mathbf{X}_{\mathbf{j}_{(a)}}) \right] \tag{18}$$

*and*

$$\pi[T, SH_N](\mathbf{X}_{\mathbf{j}_{(N)}}) = \sum_{a=1}^N \sum_{\mathbf{j}_{(a)} \subset \mathbf{j}_{(N)}} \theta_N^{(N,a)} [T]_{N,a}^{(a)}(\mathbf{X}_{\mathbf{j}_{(a)}}) = \phi_T^{(N)}(\mathbf{X}_{\mathbf{j}_{(N)}}).$$

*Moreover, for every* $s$, $[\phi_T^{(s)}]_{s,s-1}^{(r)}(\mathbf{X}_{\mathbf{j}_{(s-1)}}) = 0$, *a.s.-*$\mathbb{P}$, *for any* $\mathbf{j}_{(s-1)} \in V_\infty(s-1)$ *and any* $0 \leq r \leq s - 1$.

Note that such a result can be applied to non-centered $T \in L_s^2(\mathbf{X}_{\mathbf{j}_{(N)}})$ by considering the statistic $T' = T - \mathbb{E}(T)$. The last relation in Proposition 1 implies that the sequence $\mathbf{X}$ is *weakly independent*, as defined in [18], Section 4. We now note a property of the coefficients $\theta^{(\cdot,\cdot)}$ that will be very useful in the sequel and which can be proven by a standard recurrence argument.

**Lemma 1.** *For every $k = 1, \ldots, N$ and every $a = 1, \ldots, k$, there exists a real number $\theta^{(k,a)}$ such that*

$$\lim_{N \to +\infty} \binom{N}{k} \theta^{(k,a)}_{N*} = \theta^{(k,a)}. \tag{19}$$

For instance, by using the computations contained in [18], Section 6, we obtain

$$\theta^{(1,1)} = \alpha(A) + 1;$$

$$\theta^{(2,1)} = (\alpha(A) + 3)(\alpha(A) + 2); \qquad \theta^{(2,2)} = \frac{(\alpha(A) + 3)(\alpha(A) + 1)}{2}. \tag{20}$$

We stress that each of the coefficients $\theta^{(k,a)}_{N*}$ can be calculated in a finite number of steps. Below, it will be proven that the coefficients $\theta^{(k,a)}$ are those appearing in formula (10) above. To conclude, we present a characterization of symmetric Hoeffding spaces that is implicitly proved in [18].

**Proposition 2.** *Using the notation and assumptions of Proposition 1, for any $i = 1, \ldots, N$, a centered random variable $T \in L^2_s(\mathbf{X}_{\mathbf{j}_{(N)}})$ is an element of $SH_i(\mathbf{X}_{\mathbf{j}_{(N)}})$ if and only if there exists $h^{(i)}$ on $A^i$ such that (1) $h^{(i)} \in \Xi_i(\mathbf{X})$ and (2)*

$$T(\mathbf{X}_{\mathbf{j}_{(N)}}) = \sum_{\mathbf{j}_{(i)} \subset \mathbf{j}_{(N)}} h^{(i)}(\mathbf{X}_{\mathbf{j}_{(i)}}).$$

*Moreover, the function $h^{(i)}$ is unique in the sense that if $h'^{(i)}$ also satisfies conditions (1)–(2) above, then, $\mathbb{P}$-a.s.,*

$$h^{(i)}(\mathbf{X}_{\mathbf{j}_{(i)}}) = h'^{(i)}(\mathbf{X}_{\mathbf{j}_{(i)}}).$$

Note that the first part of the statement of Proposition 2 implies that the sequence $\mathbf{X}$ is *Hoeffding decomposable*, in the sense of [18], Definition 1.

## 4. Main results

### 4.1. Preliminaries on multiple integrals with respect to DF processes

Let $D$ be a DF process with parameter $\alpha$. We shall explore the properties of objects of the form

$$\int_{A^n} h^{(n)} \, \mathrm{d}D^{\otimes n}, \tag{21}$$

where $h^{(n)}$ is a real-valued function on $A^n$ such that

$$\mathbb{E}[h^{(n)}(\mathbf{X}_n)] = 0 \tag{22}$$

and

$$h^{(n)}(\mathbf{X}_n) \in L_s^2(\mathbf{X}_n), \qquad (23)$$

where, as before, $\mathbf{X}_n$ represents the first $n$ instants of the associated Pólya sequence with parameter $\alpha$.

**Remark.** Note that considering multiple integrals of symmetric functions is not a restriction. As a matter of fact, let $f^{(n)}$ be a measurable and not necessarily symmetric function on $A^n$. The symmetry of the product measure then yields that

$$\int_{A^n} f^{(n)} \, \mathrm{d}D^{\otimes n} = \int_{A^n} \widetilde{f}^{(n)} \, \mathrm{d}D^{\otimes n},$$

where

$$\widetilde{f}^{(n)}(a_1, \ldots, a_n) = \frac{1}{n!} \sum_\sigma f^{(n)}(a_{\sigma(1)}, \ldots, a_{\sigma(n)})$$

and $\sigma = (\sigma(1), \ldots, \sigma(n))$ runs over all permutations of $(1, \ldots, n)$.

Objects such as (21) define square integrable random variables whose variances and covariances are explicitly known. To see this, we use the notation of the previous section and write

$$h^{(n)}(\mathbf{X}_n) = \sum_{s=1}^n \pi[h^{(n)}, SH_s](\mathbf{X}_n) = \sum_{s=1}^n \sum_{\mathbf{j}_{(s)} \in V_n(s)} \phi_{h^{(n)}}^{(s)}(\mathbf{X}_{\mathbf{j}_{(s)}}), \qquad \mathbb{P}\text{-a.s.},$$

where $\pi[h^{(n)}, SH_s]$ is the projection on the $s$th symmetric Hoeffding space generated by $\mathbf{X}_n$ and the function $\phi_{h^{(n)}}^{(s)} \in \Xi_s(\mathbf{X})$ is given by applying formula (18) to the r.v. $h^{(n)}(\mathbf{X}_n)$, regarded as a symmetric element of $L^2(\mathbf{X}_n)$. This immediately yields the following identity, again due to the symmetry of the product measure $D^{\otimes n}$:

$$\int_{A^n} h^{(n)} \, \mathrm{d}D^{\otimes n} = \sum_{s=1}^n \binom{n}{s} \int_{A^s} \phi_{h^{(n)}}^{(s)} \, \mathrm{d}D^{\otimes s}. \qquad (24)$$

**Proposition 3 (Covariance between multiple integrals of the same order).** *For $n \geq 1$, let $f^{(n)}$ and $h^{(n)}$ be real-valued and symmetric functions on $A^n$ satisfying conditions (22) and (23). Then,*

$$\mathbb{E}\left[\int_{A^n} h^{(n)} \, \mathrm{d}D^{\otimes n} \int_{A^n} f^{(n)} \, \mathrm{d}D^{\otimes n}\right] = \sum_{s=1}^n \binom{n}{s}^2 \prod_{l=1}^s \frac{s-l+1}{\alpha(A)+s+l-1} \mathbb{E}[\phi_{h^{(n)}}^{(s)}(\mathbf{X}_s)\phi_{f^{(n)}}^{(s)}(\mathbf{X}_s)].$$

**Proof.** This is a consequence of Corollary 8 in [18]. Start by writing

$$\mathbb{E}\left[\int_{A^n} h^{(n)} \, \mathrm{d}D^{\otimes n} \int_{A^n} f^{(n)} \, \mathrm{d}D^{\otimes n}\right] = \sum_{s=1}^{n}\sum_{t=1}^{n} \binom{n}{s}\binom{n}{t} \mathbb{E}\left[\int_{A^s} \phi_{h^{(n)}}^{(t)} \, \mathrm{d}D^{\otimes t} \int_{A^s} \phi_{f^{(n)}}^{(s)} \, \mathrm{d}D^{\otimes s}\right]$$

$$= \sum_{s=1}^{n}\sum_{t=1}^{n} \binom{n}{s}\binom{n}{t} \mathbb{E}\left[\int_{A^{s+t}} \phi_{h^{(n)}}^{(t)} \phi_{f^{(n)}}^{(s)} \, \mathrm{d}D^{\otimes t+s}\right]$$

and observe that $D$ is the directing measure of $\mathbf{X}$, yielding that, for every $s$, $t$,

$$\mathbb{E}\left[\int_{A^{s+t}} \phi_{h^{(n)}}^{(t)} \phi_{f^{(n)}}^{(s)} \, \mathrm{d}D^{\otimes t+s}\right]$$

$$= \mathbb{E}[\phi_{h^{(n)}}^{(t)}(\mathbf{X}_{\mathbf{j}_{(s)}})\phi_{f^{(n)}}^{(s)}(\mathbf{X}_{\mathbf{i}_{(t)}})]$$

$$= \begin{cases} 0, & \text{if } t \neq s, \\ \prod_{l=1}^{s} \dfrac{s-l+1}{\alpha(A)+s+l-1} \mathbb{E}[\phi_{f^{(n)}}^{(s)}(\mathbf{X}_{\mathbf{j}_{(s)}})\phi_{h^{(n)}}^{(s)}(\mathbf{X}_{\mathbf{j}_{(s)}})], & \text{if } t = s, \end{cases}$$

where $\mathbf{j}_{(s)} \in V_{\infty}(s)$, $\mathbf{i}_{(t)} \in V_{\infty}(t)$ and $\mathrm{Card}(\mathbf{j}_{(s)} \wedge \mathbf{i}_{(t)}) = 0$, due to the fact that $\phi_{f^{(n)}}^{(s)}, \phi_{h^{(n)}}^{(s)} \in \Xi_s(\mathbf{X})$ for every $s$, as well as Corollary 8 in [18]. □

By choosing $f^{(n)} = h^{(n)}$, we immediately obtain that random variables such as (21) are square integrable. Note, moreover, that Proposition 3 contains, as a very special case, the classic computations of [8], Theorem 4. We now show that multiple integrals of the same order are almost uniquely defined.

**Proposition 4.** *Let $h^{(n)}$ and $f^{(n)}$ be symmetric functions on $A^n$ that satisfy (22) and (23). The condition*

$$\int_{A^n} h^{(n)} \, \mathrm{d}D^{\otimes n} = \int_{A^n} f^{(n)} \, \mathrm{d}D^{\otimes n}, \qquad \mathbb{P}\text{-}a.s.$$

*then implies that, for every $\mathbf{j}_{(n)} \in V_{\infty}(n)$,*

$$h^{(n)}(\mathbf{X}_{\mathbf{j}_{(n)}}) = f^{(n)}(\mathbf{X}_{\mathbf{j}_{(n)}}), \qquad \mathbb{P}\text{-}a.s.$$

**Proof.** We first consider the case where both $h^{(n)}$ and $f^{(n)}$ are elements of $\Xi_n(\mathbf{X})$. To prove the result in this case, we simply observe that the assumptions and Proposition 3 imply that

$$\mathbb{E}\left[\left(\int_{A^n} f^{(n)} \, \mathrm{d}D^{\otimes n}\right)^2\right] = \mathbb{E}\left[\left(\int_{A^n} h^{(n)} \, \mathrm{d}D^{\otimes n}\right)^2\right]$$

$$= \mathbb{E}\left[\left(\int_{A^n} h^{(n)} \, \mathrm{d}D^{\otimes n}\right)\left(\int_{A^n} f^{(n)} \, \mathrm{d}D^{\otimes n}\right)\right]$$

$$= c(n, \alpha(A))\mathbb{E}[h^{(n)}(\mathbf{X}_n)f^{(n)}(\mathbf{X}_n)]$$
$$= c(n, \alpha(A))\mathbb{E}[h^{(n)}(\mathbf{X}_n)^2]$$
$$= c(n, \alpha(A))\mathbb{E}[f^{(n)}(\mathbf{X}_n)^2],$$

where, for any $n \geq 1$,

$$c(n, \alpha(A)) := \prod_{l=1}^{n} \frac{n-l+1}{\alpha(A)+n+l-1} > 0, \tag{25}$$

thus yielding

$$\mathbb{E}[(h^{(n)}(\mathbf{X}_n) - f^{(n)}(\mathbf{X}_n))^2] = 0.$$

Now, given general $f^{(n)}$ and $h^{(n)}$, as in the statement, we write

$$\int_{A^n} h^{(n)} \, \mathrm{d}D^{\otimes n} = \sum_{s=1}^{n} \binom{n}{s} \int_{A^s} \phi_{h^{(n)}}^{(s)} \, \mathrm{d}D^{\otimes s} \quad \text{and} \quad \int_{A^n} f^{(n)} \, \mathrm{d}D^{\otimes n} = \sum_{s=1}^{n} \binom{n}{s} \int_{A^s} \phi_{f^{(n)}}^{(s)} \, \mathrm{d}D^{\otimes s},$$

using the notation of formula (24), so that the relation

$$\int_{A^n} h^{(n)} \, \mathrm{d}D^{\otimes n} = \int_{A^n} f^{(n)} \, \mathrm{d}D^{\otimes n}, \qquad \mathbb{P}\text{-a.s.},$$

implies that

$$0 = \mathbb{E}\left[\left(\int_{A^n} (h^{(n)} - f^{(n)}) \, \mathrm{d}D^{\otimes n}\right)^2\right] = \sum_{s=1}^{n} \binom{n}{s}^2 c(s, \alpha(A))\mathbb{E}[(\phi_{h^{(n)}}^{(s)} - \phi_{f^{(n)}}^{(s)})^2(\mathbf{X}_n)].$$

This immediately yields

$$\phi_{h^{(n)}}^{(s)}(\mathbf{X}_n) = \phi_{f^{(n)}}^{(s)}(\mathbf{X}_n), \qquad \mathbb{P}\text{-a.s.},$$

and therefore, $\mathbb{P}$-a.s.,

$$h^{(n)}(\mathbf{X}_n) = \sum_{s=1}^{n} \sum_{\mathbf{j}_{(s)} \in V_n(s)} \phi_{h^{(n)}}^{(s)}(\mathbf{X}_n) = \sum_{s=1}^{n} \sum_{\mathbf{j}_{(s)} \in V_n(s)} \phi_{f^{(n)}}^{(s)}(\mathbf{X}_n) = f^{(n)}(\mathbf{X}_n). \qquad \square$$

We now introduce the following class of subspaces of $L^2(D)$: $\mathcal{M}_0(D) = \Re$ and, for every $n \geq 1$,

$$\mathcal{M}_n(D) = \left\{ Y \in L^2(D) : Y = \int_{A^n} h^{(n)} \, \mathrm{d}D^{\otimes n} \text{ and } h^{(n)} \in \Xi_n(\mathbf{X}) \right\}. \tag{26}$$

By Proposition 3, it is immediate that, for any $n$, the set $\mathcal{M}_n(D)$ is an $L^2$-closed vector space isomorphic to $\sqrt{c(n, \alpha(A))}\Xi_n(\mathbf{X})$, which is, by definition, isomorphic to $\sqrt{c(n,\alpha(A))}SH_n(\mathbf{X}_n)$, that is, to the $n$th symmetric Hoeffding space associated with $\mathbf{X}_n$, endowed with the modified inner product $c(n,\alpha(A)) \times \langle\cdot,\cdot\rangle_{L^2(\mathbf{X})}$, where $c(n,\alpha(A))$ is defined in (25). Moreover, $\mathcal{M}_n(D) \perp \mathcal{M}_k(D)$ in $L^2(D)$ for every $k \neq n$. As a matter of fact, if we let $n > k$ and suppose that $h^{(n)} \in \Xi_n(\mathbf{X})$ and $h^{(k)} \in \Xi_k(\mathbf{X})$, then

$$\mathbb{E}\left[\int_{A^n} h^{(n)} \, \mathrm{d}D^{\otimes n} \int_{A^k} h^{(k)} \, \mathrm{d}D^{\otimes k}\right] = \mathbb{E}\left[\int_{A^{k+n}} h^{(n)} h^{(k)} \, \mathrm{d}D^{\otimes k+n}\right]$$
$$= \mathbb{E}[h^{(n)}(X_1,\ldots,X_n) h^{(k)}(X_{n+1},\ldots,X_{n+k})]$$
$$= \mathbb{E}[\mathbb{E}(h^{(n)}(X_1,\ldots,X_n) \mid X_{n+1},\ldots,X_{n+k})$$
$$\times h^{(k)}(X_{n+1},\ldots,X_{n+k})] = 0.$$

Proposition 4 ensures that every element $\mathcal{M}_n(D)$ admits a unique representation as an integral of an element of $\Xi_n(\mathbf{X})$ with respect to $D^{\otimes n}$.

**Remark.** In general, if a random variable $Z \in L^2(D)$ admits a representation of the type $Z = \int_{A^n} f^{(n)} \, \mathrm{d}D^{\otimes n}$, where $f^{(n)}$ is symmetric and satisfies (23), then $Z$ can be also written as $Z = \int_{A^{n+1}} \overline{f}^{(n+1)} \, \mathrm{d}D^{\otimes n+1}$, where

$$\overline{f}^{(n+1)}(a_1,\ldots,a_{n+1}) = \frac{1}{(n+1)!} \sum_\sigma f^{(n)}(a_{\sigma(1)},\ldots,a_{\sigma(n)}) \mathbf{1}_A(a_{\sigma(n+1)})$$

and $\sigma$ runs over all permutations of the class $\{1,\ldots,n+1\}$. However, the orthogonality relations discussed above ensure that if there exist $f^{(n)} \in \Xi_n(\mathbf{X})$ and $f^{(n+1)} \in \Xi_{n+1}(\mathbf{X})$ such that

$$Z = \int_{A^n} f^{(n)} \, \mathrm{d}D^{\otimes n} = \int_{A^{n+1}} f^{(n+1)} \, \mathrm{d}D^{\otimes n+1},$$

then, necessarily, $Z = 0$, $\mathbb{P}$-a.s., and therefore, by Proposition 4, $f^{(n)} = f^{(n+1)} = 0$. This also implies that if

$$F = \int_{A^n} h^{(n)} \, \mathrm{d}D^{\otimes n} = \int_{A^{n+1}} h^{(n+1)} \, \mathrm{d}D^{\otimes n+1},$$

where $h^{(k)}$, $k = n, n+1$, are symmetric and satisfy (22) and (23), then, necessarily, $\pi[h^{(n+1)}, SH_{n+1}](\mathbf{X}_{n+1}) = 0$, $\mathbb{P}$-a.s.

## 4.2. Chaotic decomposition of $L^2(D)$: $U$-statistics and polynomial complexity

We shall now show that the sequence $\{\mathcal{M}_n(D) : n \geq 0\}$ defines an orthogonal decomposition of the space $L^2(D)$. To see this, we state and prove a simple result.

**Lemma 2.** *On a probability space* $(S, \mathcal{S}, \mathbb{Q})$, *let* $\nu(\cdot; \omega)$ *be a random non-negative measure on* $(A, \mathcal{A})$ *and denote by* $L^2(\nu)$ *the class of square integrable functionals of* $\nu$. *Suppose that there exists* $q > 0$ *such that*

$$\sup_{C \in \mathcal{A}} \nu(C; \cdot) \leq q, \qquad \mathbb{Q}\text{-}a.s. \tag{27}$$

*Then, the class of random variables,*

$$\left\{ \prod_{j=1}^{n} (\nu(C_j))^{k_j} : n \geq 0, C_1, \ldots, C_n \in \mathcal{A}, 1 \leq k_1 \leq \cdots \leq k_n < +\infty \right\}$$

*is total in* $L^2(\nu)$.

**Proof.** We simply note that (27) implies that, for every $C_1, \ldots, C_n \in \mathcal{A}$ and every $(\lambda_1, \ldots, \lambda_n) \in \Re^n$,

$$\exp\left( \sum_{j=1}^{n} \lambda_j \nu(C_j) \right) = L^2 \text{-} \lim_{k \to +\infty} \prod_{j=1}^{n} \left[ \sum_{i=0}^{k} \frac{(\lambda_j \nu(C_j))^i}{i!} \right]$$

and that r.v.'s of the type $\exp(\sum_{i=1}^{n} \lambda_i \nu(C_i))$ are trivially total in $L^2(\nu)$. $\qquad \square$

We now come to the main result of this subsection.

**Theorem 1.** *Let $D$ be a DF process of parameter $\alpha$. Every $F \in L^2(D)$ then admits a unique representation of the type (6), with $h_{(F,n)} \in \Xi_n(\mathbf{X})$ for every $n \geq 1$.*

**Proof.** Due to Lemma 2, and since $D$ is a random probability measure, it is sufficient to show that random variables of the form

$$F = \prod_{j=1}^{n} (D(C_j))^{k_j}, \qquad C_1, \ldots, C_n \in \mathcal{A},$$

admit the representation (6). But,

$$F = \mathbb{E}(F) + \int_{A^{K_n}} f_{(C_1, \ldots, C_n)} \, \mathrm{d} D^{\otimes K_n},$$

where $K_t := k_1 + \cdots + k_t$ for $t = 0, \ldots, n$ (of course, $K_0 = 0$) and

$$f_{(C_1, \ldots, C_n)}(a_1, \ldots, a_{K_n}) := \frac{1}{K_n!} \sum_{\sigma} \prod_{j=1}^{n} \prod_{l=1+K_{j-1}}^{K_j} \mathbf{1}_{C_j}(a_{\sigma(l)}) - \mathbb{E}(F),$$

with $\sigma$ running over all permutations of $(1, \dots, K_n)$. We now apply formula (18) to the function $f_{(C_1,\dots,C_n)}$ and obtain

$$\int_{A^{K_n}} f_{(C_1,\dots,C_n)} \, \mathrm{d}D^{\otimes K_n} = \sum_{s=1}^{K_n} \binom{K_n}{s} \int_{A^s} \phi_{f_{(C_1,\dots,C_n)}}^{(s)} \, \mathrm{d}D^{\otimes s},$$

which implies that $F$ admits a decomposition such as (6), with

$$h_{(F,s)}(a_1, \dots, a_s) = \binom{K_n}{s} \phi_{f_{(C_1,\dots,C_n)}}^{(s)}(a_1, \dots, a_s)$$

for $s \leq K_n$ and $h_{(F,s)} = 0$ for $s > K_n$. The general result is achieved by using a standard density argument as well as the fact that each $\mathcal{M}_s(D)$ is an $L^2$-closed vector space.   □

**Remarks.** (a) Since $D$ is the directing measure of the sequence $\mathbf{X} = \{X_n : n \geq 1\}$, for every $n \geq 1$ and every $h_{(n)} \in L^2(\mathbf{X}_n)$, we have $\int_{A^n} h_{(n)} \, \mathrm{d}D^{\otimes n} = \mathbb{E}[h_{(n)}(\mathbf{X}_n) \mid D]$. It follows that, for every $F \in L^2(D)$, formula (6) can be rewritten as

$$F = \mathbb{E}(F) + \sum_{n \geq 1} \mathbb{E}[h_{(F,n)}(\mathbf{X}_n) \mid D]. \tag{28}$$

(b) We may obtain a representation similar to (28) by using elementary martingale theory. Indeed, consider a random variable $H \in L^2(\mathbf{X})$, as well as the filtration $\mathcal{X}_n = \sigma(\mathbf{X}_n)$, $n \geq 0$. It is clear that the process $Y_n = \mathbb{E}[H \mid \mathcal{X}_n]$ is a square-integrable $\mathcal{X}_n$-martingale such that $Y_n \overset{L^2}{\to} H$ as $n \to +\infty$. Now, define $g_{(H,n)}(\mathbf{X}_n) = \mathbb{E}[H \mid \mathcal{X}_n] - \mathbb{E}[H \mid \mathcal{X}_{n-1}]$, $n \geq 0$, so that we obtain immediately that

$$H = \mathbb{E}(H) + \sum_{n \geq 1} g_{(H,n)}(\mathbf{X}_n),$$

where the series on the right converges in $L^2(\mathbf{X})$. As a consequence, by conditioning with respect to $D$, we obtain

$$\mathbb{E}[H \mid D] - \mathbb{E}[\mathbb{E}[H \mid D]] = \sum_{n \geq 1} \mathbb{E}[g_{(H,n)}(\mathbf{X}_n) \mid D] = \sum_{n \geq 1} \int_{A^n} g_{(H,n)} \, \mathrm{d}D^{\otimes n}. \tag{29}$$

However, the above representation of $\mathbb{E}[H \mid D]$ is not "chaotic" in the proper sense. As a matter of fact, since the kernels $g_{(H,n)}$ are, in general, not symmetric, the integrals $\int_{A^n} g_{(H,n)} \, \mathrm{d}D^{\otimes n}$ appearing in (29) may be *not* orthogonal in $L^2(D)$, although $\mathbb{E}[g_{(H,n)}(\mathbf{X}_n)g_{(H,m)}(\mathbf{X}_m)] = 0$, for $m \neq n$. To see this, simply consider $H = h(\mathbf{X}_2)$ such that $\pi[H, SH_1] \neq 0$.

In what follows, we give two characterizations of the spaces $\mathcal{M}_j(D)$: in terms of polynomial complexity and in terms of $U$-statistics based on finite samples of the underlying sequence $\mathbf{X}$. Both are related to the following result.

**Lemma 3.** *Let $\mu_n$ be a symmetric and finite measure on the product space $(A^n, \mathcal{A}^n)$, $n \geq 2$, that is,*

$$\int_{A^n} [\mathbf{1}_{C_1} \times \cdots \times \mathbf{1}_{C_n}] \, \mathrm{d}\mu_n = \int_{A^n} [\mathbf{1}_{C_{\sigma(1)}} \times \cdots \times \mathbf{1}_{C_{\sigma(n)}}] \, \mathrm{d}\mu_n$$

*for every $C_1, \ldots, C_n \in \mathcal{A}$ and every permutation $\sigma$ of $(1, \ldots, n)$, and denote by $L_s^2(\mu_n)$ the space of symmetric and measurable functions $f$ on $A^n$ such that*

$$\int_{A^n} f^2 \, \mathrm{d}\mu_n < +\infty.$$

*The class*

$$\mathcal{H}_n := \left\{ h : h(a_1, \ldots, a_n) = \sum_\sigma \prod_{j=1}^n \mathbf{1}_{C_{\sigma(j)}}(a_j), C_j \in \mathcal{A}, j = 1, \ldots, n \right\}, \qquad (30)$$

*where $\sigma$ in the summation runs over all permutations of $\{1, \ldots, n\}$, is then total in $L_s^2(\mu_n)$.*

**Proof.** Consider an element $g \in L_s^2(\mu_n)$ such that $g \perp \mathcal{H}_n$. This means that

$$\sum_\sigma \int_{A^n} g \prod_{j=1}^n \mathbf{1}_{C_{\sigma(j)}} \, \mathrm{d}\mu_n = 0$$

for every $C_1, \ldots, C_n \in \mathcal{A}$. Now, the left side of the preceding formula equals

$$n! \int_{A^n} g \prod_{j=1}^n \mathbf{1}_{C_j} \, \mathrm{d}\mu_n,$$

due to the symmetry of $\mu_n$ and $g$. However, functions such as $\prod_{j=1}^n \mathbf{1}_{C_j}$ are total in $L^2(\mu_n)$, that is, the space of square-integrable (and not necessarily symmetric) functions on $A^n$. This implies that $g = 0$, $\mu_n$-a.s., and the result is therefore proved. $\qquad \square$

In what follows, for every DF process $D$, we define $\mathcal{P}_0(D) := \Re$ and, for $N \geq 1$, we set

$$\mathcal{P}_N(D) := \overline{\text{v.s.} \left\{ \prod_{j=1}^n (D(C_j))^{k_j} : n \geq 0, C_1, \ldots, C_n \in \mathcal{A}, k_j \geq 1, K_n \leq N \right\}}^{L^2(D)},$$

where $K_n = \sum_{j=1,\ldots,n} k_j$, to be the closed vector space generated by the polynomial functionals of $D$ with order less than or equal to $N$. Note that the sequence $\{\mathcal{P}_N(D) : N \geq 0\}$ generalizes the class $\{\mathcal{P}_n(\eta) : n \geq 0\}$ defined in the Introduction. In particular, $\mathcal{P}_N(D) \subset \mathcal{P}_{N+1}(D)$ for every $N \geq 0$ and, again be to Lemma 2, the union of the $\mathcal{P}_N(D)$'s is dense in $L^2(D)$. Consistent with the notation introduced above, we will also write $\mathcal{J}_0(D) := \Re$ and, for $N \geq 1$, $\mathcal{J}_N(D) := \mathcal{P}_N(D) \cap \mathcal{P}_{N-1}(D)^\perp$. The next result shows that the elements of $\mathcal{M}_n(D)$ are the analogues of Jacobi polynomials for general DF processes.

**Proposition 5.** *Using the above notation and assumptions, for every $N \geq 0$,*

$$\mathcal{J}_N(D) = \mathcal{M}_N(D) \quad and \quad \bigoplus_{l=0}^N \mathcal{M}_l(D) = \mathcal{P}_N(D). \tag{31}$$

**Proof.** As a by-product of the proof of Theorem 1, we know that every element of $\mathcal{J}_N(D)$ also belongs to $\bigoplus_{l=0}^N \mathcal{M}_l(D)$. Now, consider, for simplicity, a centered $F \in L^2(D)$ such that

$$F = \sum_{n=1}^N \int_{A^n} h_{(F,n)} \, \mathrm{d}D^{\otimes n}$$

with $h_{(F,n)} \in \Xi_n(\mathbf{X})$ and suppose that $F \perp \mathcal{P}_N(D)$. This implies, in particular, that, for every $C \in \mathcal{A}$,

$$\begin{aligned}
0 &= \mathbb{E}\left[\int_A h_{(F,1)} \, \mathrm{d}D \int_A \mathbf{1}_C \, \mathrm{d}D\right] \\
&= \mathbb{E}[h_{(F,1)}(X_1)\mathbf{1}_C(X_2)] \\
&= \mathbb{E}[h_{(F,1)}(X_1)(\mathbf{1}_C(X_2) - \mathbb{P}(X_2 \in C))] \\
&= c(1, \alpha(A))\mathbb{E}[h_{(F,1)}(X_1)(\mathbf{1}_C(X_1) - \mathbb{P}(X_1 \in C))] \\
&= c(1, \alpha(A))\mathbb{E}[h_{(F,1)}(X_1)\mathbf{1}_C(X_1)],
\end{aligned}$$

where the coefficients $c(n, \alpha(A))$, $n \geq 1$, are given by (25), thus yielding $h_{(F,1)}(X_1) = 0$, $\mathbb{P}$-a.s. We now use a recurrence argument. Suppose that we have shown that $F \perp \mathcal{P}_N(D)$ implies $h_{(F,j)} = 0$ for every $j = 1, \ldots, n-1$, where $n \leq N$. We then have, necessarily, for every $C_1, \ldots, C_n \in \mathcal{A}$, with the same notation as in the previous section,

$$0 = \mathbb{E}\left[\int_{A^n} h_{(F,n)} \, \mathrm{d}D^{\otimes n} \prod_{j=1}^n D(C_j)\right] = \mathbb{E}\left[\int_{A^n} h_{(F,n)} \, \mathrm{d}D^{\otimes n} \int_{A^n} \left(\prod_{j=1}^n \mathbf{1}_{C_j}\right)^{\sim} \mathrm{d}D^{\otimes n}\right],$$

where, with the usual notation,

$$\left(\prod_{j=1}^n \mathbf{1}_{C_j}\right)^{\sim}(a_1, \ldots, a_n) = \frac{1}{n!}\sum_\sigma \prod_{j=1}^n \mathbf{1}_{C_j}(a_{\sigma(j)})$$

and therefore

$$0 = \mathbb{E}\left[h_{(F,n)}(X_1, \ldots, X_n) \times \pi\left[\left(\prod_{j=1}^n \mathbf{1}_{C_j}\right)^{\sim}, SH_n\right](X_{n+1}, \ldots, X_{2n})\right]$$

$$= c(n, \alpha(A)) \mathbb{E}\left[ h_{(F,n)}(X_1, \ldots, X_n) \times \pi\left[ \left( \prod_{j=1}^{n} \mathbf{1}_{C_j} \right)^{\sim}, SH_n \right](X_1, \ldots, X_n) \right]$$

$$= \frac{c(n, \alpha(A))}{n!} \mathbb{E}\left[ h_{(F,n)}(X_1, \ldots, X_n) \sum_{\sigma} \prod_{j=1}^{n} \mathbf{1}_{C_j}(X_{\sigma(j)}) \right],$$

thus implying, by Lemma 3 and the fact that the law of $(X_1, \ldots, X_n)$ induces a symmetric measure on $A^n$, $h_{(F,n)}(\mathbf{X}_n) = 0$, $\mathbb{P}$-a.s. The proof is completed by means of standard arguments. $\qquad\square$

As announced, we have a second representation for the family $\{\mathcal{M}_n(D) \colon n \geq 0\}$.

**Proposition 6.** *Using the previous notation, let $D$ be a DF process with parameter $\alpha$ and consider the associated Pólya sequence with the same parameter $\alpha$. For every $N \geq 1$, the space $\bigoplus_{l=0}^{N} \mathcal{M}_l(D) = \mathcal{P}_n(D)$ is then generated by random variables with the representation*

$$Y = L^2\text{-}\lim_{K \to +\infty} \frac{1}{\binom{K}{N}} \sum_{\mathbf{j}_{(N)} \in V_K(N)} h(\mathbf{X}_{\mathbf{j}_{(N)}}), \qquad h \in \mathcal{H}_N, \tag{32}$$

*where the family $\mathcal{H}_N$ is defined as in (30).*

**Proof.** We first prove that if $Y$ satisfies (32), then $Y = \int_{A^N} h \, \mathrm{d}D^{\otimes N}$. It is sufficient to prove such a claim for $N = 2$ and the general case can be achieved by a standard recurrence argument. To see this, simply choose

$$h(a_1, a_2) = [\mathbf{1}_{C_1}(a_1) \mathbf{1}_{C_2}(a_2) + \mathbf{1}_{C_1}(a_2) \mathbf{1}_{C_2}(a_1)],$$

where $C_1, C_2 \in \mathcal{A}$, so that

$$\frac{1}{2} Y = L^2\text{-}\lim_{K \to +\infty} \frac{2}{K^2} \sum_{\mathbf{j}_{(2)} \in V_K(2)} \frac{1}{2} h(\mathbf{X}_{\mathbf{j}_{(2)}})$$

$$= L^2\text{-}\lim_{K \to +\infty} \frac{2}{K^2} \sum_{\mathbf{j}_{(2)} \in V_K(2)} \frac{1}{2} h(\mathbf{X}_{\mathbf{j}_{(2)}}) + L^2\text{-}\lim_{K \to +\infty} \frac{1}{K^2} \sum_{i=1}^{K} \mathbf{1}_{C_1}(X_i) \mathbf{1}_{C_2}(X_i),$$

where the last equality follows from Blackwell–MacQueen theorem and therefore

$$\frac{1}{2} Y = L^2\text{-}\lim_{K \to +\infty} \frac{1}{K^2} \sum_{\mathbf{j}_{(2)} \in V_K(2)} [\mathbf{1}_{C_1}(X_{j_1}) \mathbf{1}_{C_2}(X_{j_2}) + \mathbf{1}_{C_1}(X_{j_2}) \mathbf{1}_{C_2}(X_{j_1})]$$

$$+ L^2\text{-}\lim_{K \to +\infty} \frac{1}{K^2} \sum_{i=1}^{K} \mathbf{1}_{C_1}(X_i) \mathbf{1}_{C_2}(X_i)$$

$$= L^2\text{-}\lim_{K \to +\infty} \frac{1}{K^2}\left(\sum_{1 \le j_1 \ne j_2 \le K} \mathbf{1}_{C_1}(X_{j_1})\mathbf{1}_{C_2}(X_{j_2}) + \sum_{i=1}^{K} \mathbf{1}_{C_1}(X_i)\mathbf{1}_{C_2}(X_i)\right)$$

$$= L^2\text{-}\lim_{K \to +\infty} \frac{1}{K^2}\left(\sum_{i=1}^{K} \mathbf{1}_{C_2}(X_i)\sum_{i=1}^{K} \mathbf{1}_{C_1}(X_i)\right)$$

$$= \int_{A^2} \mathbf{1}_{C_1}\mathbf{1}_{C_2}\, \mathrm{d}D^{\otimes 2} = \frac{1}{2}\int_{A^2} h\, \mathrm{d}D^{\otimes 2}.$$

To complete the proof, we simply observe that an application of Lemma 3, similar to the one performed in the proof of Proposition 5, implies that random variables of the type

$$\int_{A^N} h\, \mathrm{d}D^{\otimes N}, \qquad h \in \mathcal{H}_N$$

are total in $\bigoplus_{l=0}^{N} M_l$ for every $N \ge 1$.                                            $\square$

The following consequence of Proposition 6 will play an important role in the sequel.

**Corollary 2.** *For every $N \ge 1$, the space $\mathcal{M}_N(D)$ is generated by random variables of the type*

$$Y = L^2\text{-}\lim_{K \to +\infty} \frac{1}{\binom{K}{N}} \sum_{\mathbf{j}_{(N)} \in V_K(N)} h(\mathbf{X}_{\mathbf{j}_{(N)}}), \qquad h \in \Xi_N(\mathbf{X}). \qquad (33)$$

**Proof.** We use formula (18) to write, for a given $h \in \mathcal{H}_N$, the following identities:

$$F = \int_{A^N} h\, \mathrm{d}D^{\otimes N} = \mathbb{E}(F) + \sum_{i=1}^{N} \binom{N}{i} \int_{A^i} \phi_h^{(i)}\, \mathrm{d}D^{\otimes i}$$

and

$$F_K = \frac{1}{\binom{K}{N}} \sum_{\mathbf{j}_{(N)} \in V_K(N)} h(\mathbf{X}_{\mathbf{j}_{(N)}}) = \mathbb{E}(F_K) + \sum_{i=1}^{N} \frac{\binom{K-i}{N-i}}{\binom{K}{N}} \sum_{\mathbf{j}_{(i)} \in V_K(i)} \phi_h^{(i)}(\mathbf{X}_{\mathbf{j}_{(i)}}).$$

It is now clear that, since such r.v.'s as

$$\int_{A^N} h\, \mathrm{d}D^{\otimes N}, \qquad h \in \mathcal{H}_N,$$

are dense in $\bigoplus_{l=0}^{N} \mathcal{M}_l(D)$, for every $i \le N$, objects of the form

$$\int_{A^i} \phi_h^{(i)}\, \mathrm{d}D^{\otimes i}, \qquad h \in \mathcal{H}_N,$$

where the $\phi_h^{(i)}$ are defined via (18), are dense in $\mathcal{M}_i(D)$. To demonstrate the result, it is therefore sufficient to prove the following relation: for every $i \geq 1$, if $g \in \Xi_i(\mathbf{X})$, then, for every $\mathbf{j}_{(i-1)} \in V_\infty(i-1)$,

$$\mathbb{E}\left[\int_{A^i} g \, \mathrm{d}D^{\otimes i} \mid \mathbf{X}_{\mathbf{j}_{(i-1)}}\right] = 0, \qquad \mathbb{P}\text{-a.s.} \tag{34}$$

As a matter of fact, this would immediately imply that

$$\binom{N}{i} \int_{A^i} \phi_h^{(i)} \, \mathrm{d}D^{\otimes i} = L^2\text{-}\lim_{K \to +\infty} \frac{\binom{K-i}{N-i}}{\binom{K}{N}} \sum_{\mathbf{j}_{(i)} \in V_K(i)} \phi_h^{(i)}(\mathbf{X}_{\mathbf{j}_{(i)}})$$

$$= L^2\text{-}\lim_{K \to +\infty} \frac{\binom{N}{i}}{\binom{K}{i}} \sum_{\mathbf{j}_{(i)} \in V_K(i)} \phi_h^{(i)}(\mathbf{X}_{\mathbf{j}_{(i)}}).$$

To prove (34), we simply use the identity

$$\mathbb{E}\left[\int_{A^i} g \, \mathrm{d}D^{\otimes i} \mid X_{j_1} = x_1, \ldots, X_{j_{i-1}} = x_{i-1}\right] = \mathbb{E}\left[\int_{A^i} g \, \mathrm{d}\widetilde{D}^{\otimes i}\right],$$

where $\widetilde{D}$ is a DF process with parameter $\alpha + \sum_{j=1,\ldots,i-1} \delta_{x_j}$ and also

$$\mathbb{E}\left[\int_{A^i} g \, \mathrm{d}\widetilde{D}^{\otimes i}\right] = \mathbb{E}[g(X_i, \ldots, X_{2i-1}) \mid X_1 = x_1, \ldots, X_{i-1} = x_{i-1}] = 0,$$

where the first equality comes from the fact that, under the probability $\mathbb{P}[\cdot \mid X_1 = x_1, \ldots, X_{i-1} = x_{i-1}]$, the sequence $\{X_{i-1+k} : k \geq 1\}$ is a Pólya sequence with parameter $\alpha + \sum_{j=1,\ldots,i-1} \delta_{x_j}$ and the second follows from the fact that $g \in \Xi_i(\mathbf{X})$. $\qquad \square$

**Remark.** A random variable such as $Z = \binom{K}{N}^{-1} \sum_{\mathbf{j}_{(N)} \in V_K(N)} h(\mathbf{X}_{\mathbf{j}_{(N)}})$, where $h$ is as in (33), is said to be a *U-statistic with (completely) degenerate kernel $h$ of degree $N$*, based on $K$ observations. The kernel $h$ is called "degenerate" since it satisfies the condition $\mathbb{E}[h(\mathbf{X}_N) \mid \mathbf{X}_{N-1}] = 0$, $\mathbb{P}$-a.s.

### 4.3. Explicit formulae

Now, consider the coefficients $\theta^{(k,a)}$, $k \geq 1$, $1 \leq a \leq k$, that are defined in formula (19). The following result provides a way to project functionals of $D$ onto spaces of multiple integrals.

**Theorem 2.** *Using the notation and assumptions of the present section, suppose that $F \in L^2(D)$ admits the decomposition*

$$F = \mathbb{E}(F) + \sum_{n \geq 1} \int_{A^n} h_{(F,n)} \, \mathrm{d}D^{\otimes n},$$

*where $h_{(F,n)} \in \Xi_n(\mathbf{X})$. For every $n \geq 1$, $h_{(F,n)}$ must then satisfy equation (10) outside a set of measure zero with respect to the probability induced on $A^n$ by the vector $\mathbf{X}_n$.*

**Proof.** For simplicity, we consider $F$ such that $\mathbb{E}(F) = 0$. Moreover, by density, it is sufficient to demonstrate the result when

$$F = \int_{A^N} \phi_h^{(N)} \, \mathrm{d}D^{\otimes N},$$

$h \in \mathcal{H}_N$ and $\phi_h^{(N)} \in \Xi_N(\mathbf{X})$ is given by formula (18). In this case, $h_{(F,n)} = 0$ for $n \neq N$ and $h_{(F,N)} = \phi_h^{(N)}$. For $n < N$, it is true that

$$0 = \sum_{k=1}^n \theta^{(n,k)} \sum_{1 \leq j_1 < \cdots < j_k \leq n} \mathbb{E}(F \mid X_{j_1}, \ldots, X_{j_k}), \qquad \text{a.s.-}\mathbb{P},$$

since we know from the discussion contained in the proof of Corollary 2 that

$$\mathbb{E}\left(\int_{A^N} g \, \mathrm{d}D^{\otimes N} \mid X_1, \ldots, X_k\right) = 0, \qquad \text{a.s.-}\mathbb{P},$$

whenever $g \in \Xi_N(\mathbf{X})$ and $N > k$. When $n = N$, we shall use the relation

$$\int_{A^N} \phi_h^{(N)} \, \mathrm{d}D^{\otimes N} = L^2 \text{-} \lim_{K \to +\infty} \frac{1}{\binom{K}{N}} \sum_{\mathbf{j}_{(N)} \in V_K(N)} \phi_h^{(N)}(\mathbf{X}_{\mathbf{j}_{(N)}})$$

$$= L^2 \text{-} \lim_{K \to +\infty} F_K(\mathbf{X}_K)$$

that was implicitly proven in Corollary 2. Now, $F_K(\mathbf{X}_K) \in SH_N(\mathbf{X}_K)$, and we may apply formula (18) to obtain, for every $\mathbf{j}_N \in V_K(N)$,

$$\frac{1}{\binom{K}{N}} \phi_h^{(N)}(\mathbf{X}_{\mathbf{j}_{(N)}}) = \sum_{a=1}^N \theta_{K*}^{(N,a)} \sum_{\mathbf{j}_{(a)} \subset \mathbf{j}_{(N)}} \mathbb{E}(F_K(\mathbf{X}_K) \mid \mathbf{X}_{\mathbf{j}_{(a)}})$$

and therefore

$$\phi_h^{(N)}(\mathbf{X}_{\mathbf{j}_{(N)}}) = \sum_{a=1}^N \theta_{K*}^{(N,a)} \binom{K}{N} \sum_{\mathbf{j}_{(a)} \subset \mathbf{j}_{(N)}} \mathbb{E}(F_K(\mathbf{X}_K) \mid \mathbf{X}_{\mathbf{j}_{(a)}})$$

so that the conclusion is obtained in this case by letting $K$ tend to infinity and using Lemma 1, as well as the relation

$$\mathbb{E}(F \mid \mathbf{X}_{\mathbf{j}_{(a)}}) = L^2 - \lim_{K \to +\infty} \mathbb{E}(F_K(\mathbf{X}_K) \mid \mathbf{X}_{\mathbf{j}_{(a)}}).$$

We are left with the case $n > N$. Here, we shall again use the fact that, for every $K$, $F_K(\mathbf{X}_K) \in SH_N(\mathbf{X}_K)$ so that

$$0 = \sum_{a=1}^{n} \theta^{(n,a)} \sum_{\mathbf{j}_{(a)} \subset V_n(a)} \mathbb{E}(F \mid \mathbf{X}_{\mathbf{j}_{(a)}})$$

$$= L^2 - \lim_{K \to \infty} \binom{K}{n} \sum_{a=1}^{n} \theta^{(n,a)}_{K*} \sum_{\mathbf{j}_{(a)} \subset \mathbf{j}_{(n)}} \mathbb{E}(F_K(\mathbf{X}_K) \mid \mathbf{X}_{\mathbf{j}_{(a)}}),$$

is immediately proved since, due to the results contained in [18], Lemma 3 and Section 5,

$$0 = \pi[F_K, SH_n](\mathbf{X}_K) = \sum_{\mathbf{j}_{(n)} \in V_K(n)} \left[ \sum_{a=1}^{n} \theta^{(n,a)}_{K*} \sum_{\mathbf{j}_{(a)} \subset \mathbf{j}_{(n)}} \mathbb{E}(F_K(\mathbf{X}_K) \mid \mathbf{X}_{\mathbf{j}_{(a)}}) \right]. \qquad \square$$

# 5. Beta random variables and Jacobi polynomials (conclusion)

As a further illustration, we apply the results of the previous sections to the special case of the simple DF process on $\{0, 1\}$ discussed in Section 1.1.

Suppose that an urn contains $\alpha_1$ red balls and $\alpha_0$ black balls, where $\alpha_1, \alpha_0 > 0$ (we could also choose $\alpha_1$ and $\alpha_0$ to be non-integer, with a straightforward interpretation). As before, an infinite sequence of extractions is performed according to the following procedure: at each step, a ball is extracted and two balls of the same color are placed in the urn before the next extraction. We define $\mathbf{X} = \{X_n : n \geq 1\}$ to be the sequence of random variables defined as

$$X_n = \begin{cases} 1, & \text{if the } n\text{th ball extracted is red,} \\ 0, & \text{otherwise.} \end{cases}$$

It is clear that $\mathbf{X}$ is, in this case, a Pólya sequence with values in $\{0, 1\}^{\infty}$ and parameter $\alpha_1 \delta_1(\cdot) + \alpha_0 \delta_0(\cdot)$, where $\delta_x$ stands for the Dirac mass at $x$. Moreover, according to the Blackwell–MacQueen Theorem, the associated DF process $D$ has the form (2), where $\eta(\omega)$ is a Beta random variable with parameters $\alpha_1$ and $\alpha_0$. In what follows, to be consistent with the notation used in the Introduction, we will identify the DF process $D$ with the random variable $\eta$ and will therefore write $\mathcal{M}_n(\eta)$ instead of $\mathcal{M}_n(D)$, $\mathcal{P}_n(\eta)$ instead of $\mathcal{P}_n(D)$ and so on.

We also recall that, for $n \geq 1$, the space $\Xi_n(\mathbf{X})$ is defined as the class of symmetric functions $\phi$ on $\{0,1\}^n$ satisfying the relation $\mathbb{E}(\phi(\mathbf{X}_n) \mid \mathbf{X}_{n-1}) = 0$, $\mathbb{P}$-a.s., and, moreover, by easily adapting (26), the sequence $\mathcal{M}_n(\eta)$ is, in this case, such that $\mathcal{M}_0(\eta) = \Re$ and

$$\mathcal{M}_n(\eta) = \left\{ Y \in L^2(\eta) : Y = \sum_{m=0}^n \binom{n}{m} \phi(\mathbf{1}^m \mathbf{0}^{n-m}) \eta^m (1-\eta)^{n-m}, \phi \in \Xi_n(\mathbf{X}) \right\}, \qquad n \geq 1,$$

where

$$\mathbf{1}^m \mathbf{0}^{n-m} := (\underbrace{1,\ldots,1}_{m \text{ times}}, \underbrace{0,\ldots,0}_{n-m \text{ times}}).$$

The following result, already announced in the Introduction, is a consequence of Proposition 5.

**Proposition 7.** *For every $n \geq 0$, $Y \in \mathcal{M}_n(\eta)$ if and only if $Y$ is a multiple of $J_n^{\alpha_1, \alpha_0}(\eta)$, where the modified Jacobi polynomial $J_n^{\alpha_1, \alpha_0}$ is defined in (3).*

We now want to represent $J_n^{\alpha_1, \alpha_0}$ in the form $\int_{\{0,1\}^n} \phi \, \mathrm{d}D^{\otimes n}$, where $D$ is defined in (2) and $\phi \in \Xi_n(\mathbf{X})$. By Proposition 4, we know that it is sufficient to write the equation, where we let $c_{n,a} = c_{n,a}(\alpha_1 + \alpha_0 - 1, \alpha_1)$, as

$$\sum_{m=0}^n \binom{n}{m} \phi(\mathbf{1}^m \mathbf{0}^{n-m}) x^m (1-x)^{n-m} = \sum_{a=0}^n c_{n,a} x^a, \qquad x \in [0,1],$$

which is satisfied if and only if $\phi$ solves the (triangular) system

$$\sum_{m=0}^a \binom{n}{m} \binom{n-m}{n-a} (-1)^{a-m} \phi(\mathbf{1}^m \mathbf{0}^{n-m}) = c_{n,a}, \qquad a = 0, \ldots, n. \tag{35}$$

Theorem 2 and Proposition 7 yield the following result, related to formula (5) in Section 1.2.

**Proposition 8.** *Using the notation and the assumptions of this section the following conditions are equivalent:*

1. *$\psi \in \Xi_n(\mathbf{X})$;*
2. *$\psi = k\phi$, where $k$ is a real constant and $\phi$ solves the system in formula (35);*
3. *there exists $Y \in L^2(\eta)$ such that*

$$\phi(\mathbf{X}_n) = \sum_{a=1}^n \theta^{(n,a)} \sum_{\mathbf{j}_{(a)} \in V_n(a)} \mathbb{E}(Y \mid \mathbf{X}_{\mathbf{j}_{(a)}}), \qquad \mathbb{P}\text{-a.s.}, \tag{36}$$

*where the coefficients $\theta^{(\cdot, \cdot)}$ are given by formula (19).*

In other words, Proposition 8 states that, in the case of a classic Pólya urn and for every $n$, the only completely degenerate $U$-statistics of order $n$ are those constructed by means of kernels that are multiples of solutions of the system in (35). Moreover, conditional expectations of functionals of $\eta$, with respect to the underlying urn sequence $\mathbf{X}$, are linked to Jacobi polynomials via formula (36).

We conclude the section by stating the following consequence of Proposition 3.

**Proposition 9.** *For every $n \geq 1$,*

$$\int_0^1 J_n^{\alpha_1, \alpha_0}(x)^2 p_{\alpha_1, \alpha_0}(x) \, dx$$

$$= \prod_{l=1}^n \left[ \frac{n-l+1}{\alpha_1 + \alpha_0 + n + l - 1} \right] \sum_{m=0}^n \binom{n}{m} \phi(\mathbf{1}^m \mathbf{0}^{n-m})^2 \int_0^1 x^m (1-x)^{n-m} p_{\alpha_1, \alpha_0}(x) \, dx,$$

*where $\phi$ is given by (35).*

# 6. Connections with other orthogonal polynomials and multiallele diffusion models

In this section, we explain how our results can be related to *Wright–Fisher diffusion processes* (or, more generally, to *Fleming–Viot processes*) of population genetics. The reader is referred to [6], Section 10, [7] and the references therein for basic terminology and results.

We fix, here and for the rest of the section, an integer $K \geq 2$. The $K$-type Wright–Fisher process (see also [10]) is defined as the diffusion taking values in the symplex

$$\Delta_K = \left\{ \zeta_{(K)} = (\zeta_1, \ldots, \zeta_K) : \zeta_i \geq 0, \sum_{i=1}^K \zeta_i = 1 \right\}$$

and with generator

$$L = \frac{1}{2} \sum_{i,j=1}^K \zeta_i (\epsilon_{ij} - \zeta_j) \frac{\partial^2}{\partial \zeta_i \, \partial \zeta_j} + \sum_{j=1}^K \left( \sum_{i=1}^K q_{ij} \zeta_i \right) \frac{\partial}{\partial \zeta_j}, \tag{37}$$

where $\epsilon_{ij} = 1$ if $j = i$ and $= 0$ otherwise, and $(q_{ij})$ is the matrix describing the mutation structure.

In what follows, we will write

$$\Delta_{K-1}^0 = \left\{ \boldsymbol{\gamma}_{(K-1)} = (\gamma_1, \ldots, \gamma_{K-1}) : \gamma_i > 0, \sum_{i=1}^{K-1} \gamma_i < 1 \right\},$$

and, for any vector $\boldsymbol{\theta} = (\theta_1, \ldots, \theta_K)$ such that $\theta_j > 0$ $(j = 1, \ldots, K)$, we denote by $D_{\boldsymbol{\theta}}$ a DF process on $\{1, \ldots, K\}$ with parameter $\alpha_{\boldsymbol{\theta}}(\cdot) = \sum_{j=1,\ldots,K} \theta_j \delta_j(\cdot)$, where $\delta_j$ is the Dirac mass at $j$. Note that, by definition, the vector

$$D_{\boldsymbol{\theta}, K-1} = (D_{\boldsymbol{\theta}}(\{1\}), \ldots, D_{\boldsymbol{\theta}}(\{K-1\})) \tag{38}$$

is a random element such that $\mathbb{P}(D_{\boldsymbol{\theta}, K-1} \in \Delta^0_{K-1}) = 1$. We also recall that the law of $D_{\boldsymbol{\theta}, K-1}$ is absolutely continuous with respect to the restriction of the Lebesgue measure to $\Delta^0_{K-1}$ and that the associated density is

$$f_{\boldsymbol{\theta}, K-1}(\boldsymbol{\gamma}_{(K-1)}) = \frac{\Gamma(\theta_1 + \cdots + \theta_K)}{\Gamma(\theta_1) \cdots \Gamma(\theta_K)} \gamma_1^{\theta_1 - 1} \cdots \gamma_{K-1}^{\theta_{K-1} - 1} \left(1 - \sum_{j=1}^{K-1} \gamma_j\right)^{\theta_{K-1} - 1}. \tag{39}$$

We write $\mathbf{W}_{(K-1)}$ to indicate the set of vectors $\mathbf{n} = (n_1, \ldots, n_{K-1})$ such that each $n_j$ is a non-negative integer and, for each $\mathbf{n} \in \mathbf{W}_{(K-1)}$, we use the customary notation $|\mathbf{n}| := \sum_{j=1,\ldots,K-1} n_j$. A system of *biorthogonal polynomials*

$$\{\eta^{(1)}_{\mathbf{n}_1}, \eta^{(2)}_{\mathbf{n}_2}\} = \{\eta^{(1)}_{\mathbf{n}_1}, \eta^{(2)}_{\mathbf{n}_2} : \mathbf{n}_1, \mathbf{n}_2 \in \mathbf{W}_{(K-1)}\}$$

with respect to the density $f_{\boldsymbol{\theta}, K-1}$ introduced in (39) is a double collection of polynomials defined on $\Delta^0_{K-1}$, indexed by the elements of $\mathbf{W}_{(K-1)}$ and such that (a) the degree of $\eta^{(i)}_{\mathbf{n}_i}$ is given by $|\mathbf{n}_i|$ $(i = 1, 2)$ and (b) for every $\mathbf{n}_1, \mathbf{n}_2 \in \mathbf{W}_{(K-1)}$,

$$\int_{\Delta^0_{K-1}} \eta^{(1)}_{\mathbf{n}_1}(\boldsymbol{\gamma}_{(K-1)}) \eta^{(2)}_{\mathbf{n}_2}(\boldsymbol{\gamma}_{(K-1)}) f_{\boldsymbol{\theta}, K-1}(\boldsymbol{\gamma}_{(K-1)}) \, \mathrm{d}\boldsymbol{\gamma}_{(K-1)} = \begin{cases} 0, & \text{if } \mathbf{n}_1 \neq \mathbf{n}_2, \\ 1, & \text{if } \mathbf{n}_1 = \mathbf{n}_2 \end{cases}$$

(see, e.g., [21] for further details). Note that, for every $\boldsymbol{\theta} = (\theta_1, \ldots, \theta_K)$ such that $\theta_j > 0$, a *complete* system of biorthogonal polynomials with respect to $f_{\boldsymbol{\theta}, K-1}$ can be obtained by properly renormalizing the double family of generalized *Appell–Jacobi polynomials* defined in [12], formulae (2.6) and (2.7).

A crucial point in the analysis of the Wright–Fisher process defined by (37) is the explicit computation of the associated transition density. This task can be hugely simplified by introducing the additional assumption

$$q_{ij} = 2^{-1}\theta_j > 0 \qquad \forall 1 \le i \neq j \le K. \tag{40}$$

As a matter of fact, in this case, the Wright–Fisher process has a unique stationary distribution given by the law of the DF process $D_{\boldsymbol{\theta}}$ on $\{1, \ldots, K\}$, where $\boldsymbol{\theta} = (\theta_1, \ldots, \theta_K)$, introduced above. In particular, according, for example, to [12], Section 3 and Section 4, and [10], Theorem 1, assumption (40) implies that the transition density of the Wright–Fisher process can be expressed in terms of a class of *kernel orthogonal polynomials* $\{Q_n(\cdot; \cdot) : n \ge 0\}$ which is uniquely determined by the following conditions: (i) $Q_0 = 1$;

(ii) for every polynomial $R_n(\cdot)$ of degree $n \geq 1$ and defined on $\Delta_{K-1}^0$,

$$R_n(\boldsymbol{\gamma}_{(K-1)}) = \sum_{j=0}^{n} \mathbb{E}[Q_n(D_{\theta,K-1}; \boldsymbol{\gamma}_{(K-1)}) R_n(D_{\theta,K-1})] \qquad \forall \boldsymbol{\gamma}_{(K-1)} \in \Delta_{K-1}^0,$$

where $D_{\theta,K-1}$ is defined as in (38); (iii) for any complete set of biorthogonal polynomials $\{\eta_{\mathbf{m}_1}^{(1)}, \eta_{\mathbf{m}_2}^{(2)}\}$ with respect to $f_{\boldsymbol{\theta},K-1}$,

$$Q_n(\boldsymbol{\gamma}_{(K-1)}; \boldsymbol{\gamma}'_{(K-1)}) = \sum_{\mathbf{n} \in \mathbf{W}_{(K-1)} : |\mathbf{n}|=n} \eta_{\mathbf{n}}^{(1)}(\boldsymbol{\gamma}_{(K-1)}) \eta_{\mathbf{n}}^{(2)}(\boldsymbol{\gamma}'_{(K-1)}) \qquad (41)$$

for every $n \geq 1$ and every $\boldsymbol{\gamma}_{(K-1)}, \boldsymbol{\gamma}'_{(K-1)} \in \Delta_{K-1}^0$.

As already pointed out, a complete biorthogonal system of polynomials such as the one appearing in condition (iii) is explicitly computed (up to normalization constants) in [12] by means of a generalization of Appell–Jacobi polynomials. Another approach for computing the kernel polynomials $Q_n$ is used in [10], where the author uses the representation, valid for $n \geq 1$ and $\boldsymbol{\gamma}_{(K-1)}, \boldsymbol{\gamma}'_{(K-1)} \in \Delta_{K-1}^0$,

$$Q_n(\boldsymbol{\gamma}_{(K-1)}; \boldsymbol{\gamma}'_{(K-1)}) = \sum_{\mathbf{n} : |\mathbf{n}|=n} P_{\mathbf{n}}(\boldsymbol{\gamma}_{(K-1)}) P_{\mathbf{n}}(\boldsymbol{\gamma}'_{(K-1)}), \qquad (42)$$

where the family of orthogonal polynomials $\{P_{\mathbf{n}} : \mathbf{n} \in \mathbf{W}_{(K-1)}\}$ is obtained through a Gram–Schmidt orthogonalization, with respect to $f_{\boldsymbol{\theta},K-1}$, of the monomials

$$M_{\mathbf{n}}(\boldsymbol{\gamma}_{(K-1)}) = \gamma_1^{n_1} \gamma_2^{n_2} \cdots \gamma_{K-1}^{n_{K-1}}, \qquad \mathbf{n} = (n_1, \ldots, n_{K-1}) \in \mathbf{W}_{(K-1)},$$

realized by means of a total ordering of the elements of $\mathbf{W}_{(K-1)}$ (see [10], pages 311–315 for details). In particular, the degree of each $P_{\mathbf{n}}$ is equal to $|\mathbf{n}|$ and, for $n \geq 1$, the set $\{P_{\mathbf{n}} : |\mathbf{n}| \leq n\}$ is an orthonormal basis (with respect to the measure induced on $\Delta_{K-1}^0$ by the density $f_{\boldsymbol{\theta},K-1}$) of the space of polynomials of degree less than or equal to $n$.

Another consequence of our results are further probabilistic characterizations of the orthogonal families of polynomials $\{\eta_{\mathbf{m}_1}^{(1)}, \eta_{\mathbf{m}_2}^{(2)}\}$ and $\{P_{\mathbf{n}}\}$ appearing, respectively, in (41) and (42). In particular, we have the following proposition, whose proof is a standard consequence of Proposition 5 and the definition of the family $\{P_{\mathbf{n}}\}$.

**Proposition 10.** *Using the assumptions of this section and the same notation as in Proposition 5, for every $n \geq 1$ and every $F \in L^2(D_{\boldsymbol{\theta}})$, the following assertions are equivalent:*

1. *$F \in \mathcal{J}_n(D_{\boldsymbol{\theta}}) = \mathcal{M}_n(D_{\boldsymbol{\theta}})$;*
2. *there exist real constants $\{c_{\mathbf{n}} : \mathbf{n} \in \mathbf{W}_{(K-1)}, |\mathbf{n}| = n\}$ such that*

$$F = \sum c_{\mathbf{n}} P_{\mathbf{n}}(D_{\boldsymbol{\theta},K-1});$$

*in particular, the set $\{P_{\mathbf{n}}(D_{\boldsymbol{\theta},K-1}) : \mathbf{n} \in \mathbf{W}_{(K-1)}, |\mathbf{n}| = n\}$ is an orthonormal basis of $\mathcal{M}_n(D_{\boldsymbol{\theta}})$.*

Now, denote by $\mathbf{X}_n^{\boldsymbol{\theta}}$ the first $n$ instants of the Pólya sequence associated with $D_{\boldsymbol{\theta}}$ (see Section 2). It follows from Proposition 10 that there exists an orthonormal basis

$$\{h_{\mathbf{n},\boldsymbol{\theta}} : \mathbf{n} \in \mathbf{W}_{(K-1)}, |\mathbf{n}| = n\}$$

of the space $\sqrt{c(n, |\boldsymbol{\theta}|)}\,\Xi_n(\mathbf{X}_n^{\boldsymbol{\theta}})$, where $|\boldsymbol{\theta}| := \sum_{j=1}^K \theta_j$, such that, for each $\mathbf{n}$,

$$P_{\mathbf{n}}(D_{\boldsymbol{\theta},K-1}) = \int_{\{1,\ldots,K\}^n} h_{\mathbf{n},\boldsymbol{\theta}}\,\mathrm{d}D_{\boldsymbol{\theta}}^{\otimes n}, \qquad \text{a.s.-}\mathbb{P}.$$

For a fixed $\boldsymbol{\gamma}_{(K-1)} = (\gamma_1, \ldots, \gamma_{K-1}) \in \Delta_{K-1}^0$, we define the measure $\mu_{\boldsymbol{\gamma}_{(K-1)}}$ on $\{1, \ldots, K\}$ as follows:

$$\mu_{\boldsymbol{\gamma}_{(K-1)}}(\{i\}) = \gamma_i, \qquad i = 1, \ldots, K-1,$$

$$\mu_{\boldsymbol{\gamma}_{(K-1)}}(\{K\}) = 1 - \sum_{i=1}^{K-1} \mu_{\boldsymbol{\gamma}_{(K-1)}}(\{i\}).$$

For $n \geq 2$, we write $\mu_{\boldsymbol{\gamma}_{(K-1)}}^{\otimes n}$ to indicate the canonical product measure induced by $\mu_{\boldsymbol{\gamma}_{(K-1)}}$ on $\{1, \ldots, K\}^n$. Also, $\mu_{\boldsymbol{\gamma}_{(K-1)}}^{\otimes 1} := \mu_{\boldsymbol{\gamma}_{(K-1)}}$. Since

$$\mathrm{d}D_{\boldsymbol{\theta}}^{\otimes n} = \mathrm{d}\mu_{D_{\theta,K-1}}^{\otimes n}, \qquad \text{a.s.-}\mathbb{P},$$

the following characterization of the kernels $Q_n$ defined above is now easily proved.

**Proposition 11.** *Let $\{\eta_{\mathbf{n}_1}^{(1)}, \eta_{\mathbf{n}_2}^{(2)}\}$ be a complete system of biorthogonal polynomials as defined above and, for $n \geq 0$, let $Q_n(\cdot; \cdot)$ be the kernel orthogonal polynomial defined by means of conditions* (i)–(iii) *above. Then, for every $(\boldsymbol{\gamma}_{(K-1)}, \boldsymbol{\gamma}'_{(K-1)}) \in (\Delta_{K-1}^0)^2$ outside a set of zero Lebesgue measure,*

$$Q_n(\boldsymbol{\gamma}_{(K-1)}, \boldsymbol{\gamma}'_{(K-1)}) = \sum_{\mathbf{n} \in \mathbf{W}_{(K-1)}:\, |\mathbf{n}|=n} \eta_{\mathbf{n}}^{(1)}(\boldsymbol{\gamma}_{(K-1)}) \eta_{\mathbf{n}}^{(2)}(\boldsymbol{\gamma}'_{(K-1)})$$

$$= \sum_{\mathbf{n} \in \mathbf{W}_{(K-1)}:\, |\mathbf{n}|=n} \int_{\{1,\ldots,K\}^n} h_{\mathbf{n},\boldsymbol{\theta}}\,\mathrm{d}\mu_{\boldsymbol{\gamma}_{(K-1)}}^{\otimes n} \times \int_{\{1,\ldots,K\}^n} h_{\mathbf{n},\boldsymbol{\theta}}\,\mathrm{d}\mu_{\boldsymbol{\gamma}'_{(K-1)}}^{\otimes n}.$$

One can now use, for example, [10], formula (3.1), page 316, to derive an expression for the transition density of a Wright–Fisher process with generator (37) in terms of the kernels $h_{\mathbf{n},\boldsymbol{\theta}}$. Namely, we have the following.

**Corollary 3.** *The transition density at time $t > 0$ of the first $K - 1$ allele frequencies of the $K$-type Wright–Fisher diffusion process with generator (37) satisfying (40),* *given a vector of initial frequencies* $\boldsymbol{\gamma}'_{(K)} = (\gamma'_1, \ldots, \gamma'_K) \in \Delta_K$ *such that* $\boldsymbol{\gamma}'_{(K-1)} = (\gamma'_1, \ldots, \gamma'_{K-1}) \in \Delta^0_{K-1}$, *is*

$$
P(\boldsymbol{\gamma}_{(K-1)}; t, \boldsymbol{\gamma}'_{(K)})
$$

$$
= \frac{\Gamma(\theta_1 + \cdots + \theta_K)}{\Gamma(\theta_1) \cdots \Gamma(\theta_K)} \gamma_1^{\theta_1 - 1} \cdots \gamma_{K-1}^{\theta_{K-1} - 1} \left\{ 1 + \sum_{n=1}^{\infty} \rho_n(t) Q_n(\boldsymbol{\gamma}_{(K-1)}, \boldsymbol{\gamma}'_{(K-1)}) \right\}
$$

$$
= \frac{\Gamma(\theta_1 + \cdots + \theta_K)}{\Gamma(\theta_1) \cdots \Gamma(\theta_K)} \gamma_1^{\theta_1 - 1} \cdots \gamma_{K-1}^{\theta_{K-1} - 1}
$$

$$
\times \left\{ 1 + \sum_{n=1}^{\infty} \rho_n(t) \sum_{\mathbf{n} \in \mathbf{W}_{(K-1)} : |\mathbf{n}| = n} \int_{\{1, \ldots, K\}^n} h_{\mathbf{n}, \boldsymbol{\theta}} \, d\mu^{\otimes n}_{\boldsymbol{\gamma}_{(K-1)}} \times \int_{\{1, \ldots, K\}^n} h_{\mathbf{n}, \boldsymbol{\theta}} \, d\mu^{\otimes n}_{\boldsymbol{\gamma}'_{(K-1)}} \right\}
$$

*for a.e.* $\boldsymbol{\gamma}_{(K-1)} \in \Delta^0_{K-1}$, *where* $\rho_n(t) := \exp\{-\frac{1}{2}n(n-1)t - \frac{1}{2}(\theta_1 + \cdots + \theta_K)nt\}$. *In particular, if* $D_{\boldsymbol{\theta}}$ *is a DF process with parameter* $\alpha_{\boldsymbol{\theta}}$, *then for a.e.* $\boldsymbol{\gamma}_{(K-1)} \in \Delta^0_{K-1}$, *the random variable*

$$
G_{\boldsymbol{\gamma}_{(K-1)}} = P(\boldsymbol{\gamma}_{(K-1)}; t, (D_{\boldsymbol{\theta}}(\{1\}), \ldots, D_{\boldsymbol{\theta}}(\{K\})))
$$

*is an element of* $L^2(D_{\boldsymbol{\theta}})$ *and, for every* $n \geq 1$, *the projection of* $G_{\boldsymbol{\gamma}_{(K-1)}}$ *onto* $\mathcal{M}_n(D_{\boldsymbol{\theta}})$ *is given by*

$$
\frac{\Gamma(\theta_1 + \cdots + \theta_K)}{\Gamma(\theta_1) \cdots \Gamma(\theta_K)} \gamma_1^{\theta_1 - 1} \cdots \gamma_{K-1}^{\theta_{K-1} - 1}
$$

$$
\times \rho_n(t) \int_{\{1, \ldots, K\}^n} \left( \sum_{\mathbf{n} \in \mathbf{W}_{(K-1)} : |\mathbf{n}| = n} P_{\mathbf{n}}(\boldsymbol{\gamma}_{(K-1)}) \times h_{\mathbf{n}, \boldsymbol{\theta}} \right) dD^{\otimes n}_{\boldsymbol{\theta}}. \tag{43}
$$

We can finally combine (10) and (43) to obtain that, for every $(a_1, \ldots, a_n) \in \{1, \ldots, K\}^n$ and every $n \geq 1$,

$$
\frac{\Gamma(\theta_1 + \cdots + \theta_K)}{\Gamma(\theta_1) \cdots \Gamma(\theta_K)} \gamma_1^{\theta_1 - 1} \cdots \gamma_{K-1}^{\theta_{K-1} - 1} \sum_{\mathbf{n} \in \mathbf{W}_{(K-1)} : |\mathbf{n}| = n} P_{\mathbf{n}}(\boldsymbol{\gamma}_{(K-1)}) \times h_{\mathbf{n}, \boldsymbol{\theta}}(a_1, \ldots, a_n)
$$

$$
= \rho_n(t)^{-1} \sum_{k=1}^{n} \theta^{(n,k)} \sum_{1 \leq j_1 < \cdots < j_k \leq n} \mathbb{E}(G_{\boldsymbol{\gamma}_{(K-1)}} - \mathbb{E}(G_{\boldsymbol{\gamma}_{(K-1)}}) \mid X_1^{\boldsymbol{\theta}} = a_{j_1}, \ldots, X_k^{\boldsymbol{\theta}} = a_{j_k}),
$$

*where* $(X_1^{\boldsymbol{\theta}}, \ldots, X_k^{\boldsymbol{\theta}}, \ldots)$ *is the Pólya sequence associated with* $D_{\boldsymbol{\theta}}$.

To conclude, we recall that, in [7], the authors derive an explicit expression of the transition density of the Fleming–Viot process (i.e., a measure-valued generalization of the Wright–Fisher diffusion) under conditions ensuring that its stationary distribution is the law of a general DF process on a compact metric space. The relation between such a result and the orthogonal decomposition of $L^2(D)$ introduced in our paper is far from straightforward and will be investigated elsewhere.

# 7. Further examples and applications

(a) *Exponential functionals.* Consider a DF process $D$ with parameter $\alpha$. We want to write the decomposition of the functional

$$G = \exp(\lambda D(C)),$$

where $\lambda$ is a real constant and $C$ is an element of $\mathcal{A}$ such that $\alpha(A) > a(C) > 0$. This implies, according to [8], Proposition 1, that $D(C) \in (0, 1)$ with probability 1. Moreover, we know that, under $\mathbb{P}$, $D(C)$ has a Beta distribution with parameters $(\alpha(C), \alpha(A \backslash C))$. Now, the decomposition of $G$ is given by

$$G = \mathbb{E}(G) + \sum_{n \geq 1} \int_{A^n} h_{(G,n)} \, dD^{\otimes n} = {}_1F_1(\alpha(C), \alpha(A), \lambda) + \sum_{n \geq 1} \int_{A^n} h_{(G,n)} \, dD^{\otimes n},$$

where ${}_1F_1$ indicates a *confluent hypergeometric function of the first kind* (see, e.g., [1]) and

$$
\begin{aligned}
h_{(G,n)}(a_1, \ldots, a_n) &= \sum_{k=1}^{n} \theta^{(n,k)} \sum_{1 \leq j_1 < \cdots < j_k \leq n} \mathbb{E}(G - \mathbb{E}(G) \mid X_1 = a_{j_1}, \ldots, X_k = a_{j_k}) \\
&= \sum_{k=1}^{n} \theta^{(n,k)} \sum_{1 \leq j_1 < \cdots < j_k \leq n} \Bigg[ {}_1F_1\left(\alpha(C) + \sum_{i=1}^{k} \mathbf{1}_C(a_{j_i}), \alpha(A) + k, \lambda\right) \\
&\qquad\qquad\qquad\qquad\qquad\qquad - {}_1F_1(\alpha(C), \alpha(A), \lambda) \Bigg],
\end{aligned}
$$

where the second equality derives from the fact that, conditioned on $(X_{j_1}, \ldots, X_{j_k})$, $D$ is a DF process with parameter $\alpha + \sum_{i=1,\ldots,k} \delta_{X_{j_i}}$. Similar calculations apply to functionals of the type

$$G = \exp\left(\sum_{i=1,\ldots,n} \lambda_i D(C_i)\right),$$

where $(C_1, \ldots, C_n)$ is a finite measurable partition of $A$.

(b) *Approximation by U-statistics.* Now, let $D$ be a DF process with parameter $\alpha$ and let $\mathbf{X}$ be the associated Pólya sequence with parameter $\alpha$. As announced in the Introduction, we shall use the results of the previous sections to solve the problem of finding the best $L^2$ approximation of an element of $L^2(D)$ by means of a symmetric statistic of the finite sequence $\mathbf{X}_N = (X_1, \ldots, X_N)$, where $N$ is a fixed integer strictly greater than 1. Now, recall that the space of symmetric elements of $L^2(\mathbf{X}_N)$ is the direct sum of $\Re$ and the first $N$ symmetric Hoeffding spaces associated with $\mathbf{X}_N$. According to Proposition 2, this yields that the above-stated problem reduces to the following: find functions $h_i \in \Xi_i(\mathbf{X})$, $i = 1, \ldots, N$, such that

$$
\mathbb{E}\left[F - \left(\mathbb{E}(F) + \sum_{i=1}^{N} \sum_{\mathbf{j}_{(i)} \subset V_N(i)} h_i(\mathbf{X}_{\mathbf{j}_{(i)}})\right)\right]^2
$$
$$
= \operatorname*{arg\,min}_{g_i \in \Xi_i(\mathbf{X}), i=1,\ldots,N} \mathbb{E}\left[F - \left(\mathbb{E}(F) + \sum_{i=1}^{N} \sum_{\mathbf{j}_{(i)} \subset V_N(i)} g_i(\mathbf{X}_{\mathbf{j}_{(i)}})\right)\right]^2. \tag{44}
$$

To carry out this program, we introduce the coefficients, appearing in the statement of Corollary 9 in [18] and defined for $n \geq 1$ and $r = 0, \ldots, n$,

$$
c(r, n, \alpha(A)) = \prod_{l=1}^{n} \frac{n - r - l + 1}{\alpha(A) + n + l - 1}
$$

and note that

$$
c(0, n, \alpha(A)) = c(n, \alpha(A)),
$$

where the term on the right is defined in (25).

**Proposition 12.** *Suppose that $F \in L^2(D)$ admits the decomposition*

$$
F = \mathbb{E}(F) + \sum_{n \geq 1} \int_{A^n} h_{(F,n)} \, \mathrm{d}D^{\otimes n},
$$

*where $h_{(F,n)} \in \Xi_n(\mathbf{X})$. Condition (44) is then satisfied by*

$$
h_i = \frac{1}{\binom{N}{i}} h_{(F,i)}
$$

*for $i = 1, \ldots, N$. Moreover, the corresponding quadratic error is given by*

$$
\mathbb{E}\left[F - \left(\mathbb{E}(F) + \sum_{i=1}^{N} \sum_{\mathbf{j}_{(i)} \subset V_N(i)} h_i(\mathbf{X}_{\mathbf{j}_{(i)}})\right)\right]^2
$$
$$
= \sum_{n \geq N+1} c(n, \alpha(A)) \mathbb{E}[h_{(F,n)}(\mathbf{X}_n)^2]
$$
$$
+ \sum_{n=1}^{N} \mathbb{E}[h_{(F,n)}(\mathbf{X}_n)^2] \left[c(n, \alpha(A)) - \binom{N}{n}^{-1} \sum_{r=0}^{n} \binom{n}{r} \binom{N-n}{n-r}_* c(r, n, \alpha(A))\right].
$$

**Proof.** It is sufficient to prove that, for every $i \geq 1$ and every $N \geq i$, for every $h_i, g_i \in \Xi_i(\mathbf{X})$,

$$\mathbb{E}\left(\int_{A^i} h_i \, dD^{\otimes i} \sum_{\mathbf{j}_{(i)} \in V_N(i)} g_i(\mathbf{X}_{\mathbf{j}_{(i)}})\right)$$

$$= \frac{1}{\binom{N}{i}} \mathbb{E}\left(\sum_{\mathbf{j}_{(i)} \in V_N(i)} h_i(\mathbf{X}_{\mathbf{j}_{(i)}}) \sum_{\mathbf{j}_{(i)} \in V_N(i)} g_i(\mathbf{X}_{\mathbf{j}_{(i)}})\right).$$

By a density argument we can take $h_i = \phi_f^{(i)}$, given by formula (18) for a certain $f \in \mathcal{H}_i$. We may then write, using the notation $r(\mathbf{i}_{(i)}, \mathbf{j}_{(i)}) := \mathrm{Card}(\mathbf{i}_{(i)} \wedge \mathbf{j}_{(i)})$ Corollary 10 in [18], as well as a simple combinatorial argument,

$$\mathbb{E}\left(\int_{A^i} h_i \, dD^{\otimes i} \sum_{\mathbf{j}_{(i)} \in V_N(i)} g_i(\mathbf{X}_{\mathbf{j}_{(i)}})\right) = \lim_{K \to \infty} \frac{1}{\binom{K}{i}} \mathbb{E}\left(\sum_{\mathbf{j}_{(i)} \in V_K(i)} h_i(\mathbf{X}_{\mathbf{j}_{(i)}}) \sum_{\mathbf{j}_{(i)} \in V_N(i)} g_i(\mathbf{X}_{\mathbf{j}_{(i)}})\right)$$

$$= \lim_{K \to \infty} \frac{1}{\binom{K}{i}} \sum_{\mathbf{j}_{(i)} \in V_K(i)} \sum_{r=0}^{i} \sum_{\substack{\mathbf{i}_{(i)} \in V_N(i): \\ (\mathbf{j}_{(i)}, \mathbf{i}_{(i)})=r}} \mathbb{E}(h_i(\mathbf{X}_{\mathbf{j}_{(i)}}) g_i(\mathbf{X}_{\mathbf{i}_{(i)}}))$$

$$= \sum_{r=0}^{i} \binom{i}{r} \binom{N-i}{i-r}_* c(r, i, \alpha(A)) \mathbb{E}(h_i(\mathbf{X}_{\mathbf{j}_{(i)}}) g_i(\mathbf{X}_{\mathbf{j}_{(i)}}))$$

$$= \frac{1}{\binom{N}{i}} \mathbb{E}\left(\sum_{\mathbf{j}_{(i)} \in V_N(i)} h_i(\mathbf{X}_{\mathbf{j}_{(i)}}) \sum_{\mathbf{j}_{(i)} \in V_N(i)} g_i(\mathbf{X}_{\mathbf{j}_{(i)}})\right).$$

The last formula in the statement is straightforward. $\qquad\square$

For example, from the calculations in part (a), we obtain that the best approximation of $G = \exp(\lambda D(C))$, by means of $U$-statistics based on $\mathbf{X}_N$, is

$$G_N = \mathbb{E}(G) + \sum_{i=1}^{N} \sum_{\mathbf{j}_{(i)} \subset V_N(i)} h_i(\mathbf{X}_{\mathbf{j}_{(i)}})$$

$$= {}_1F_1(\alpha(C), \alpha(A), \lambda)$$

$$+ \sum_{i=1}^{N} \sum_{\mathbf{j}_{(i)} \subset V_N(i)} \sum_{k=1}^{i} \binom{N}{i}^{-1} \theta^{(i,k)}$$

$$\times \sum_{\mathbf{j}_{(k)} \subset \mathbf{j}_{(i)}} \left[ {}_1F_1\left(\alpha(C) + \sum_{l=1}^{k} \mathbf{1}_C(X_{j_l}), \alpha(A) + k, \lambda\right)\right.$$

$$- {}_1F_1(\alpha(C), \alpha(A), \lambda)\Bigg].$$

**Remark.** Suppose that $(A, \mathcal{A}) = ([0, 1], \mathcal{B}([0, 1]))$, where $\mathcal{B}$ stands for the Borel $\sigma$-field and $\alpha$ is equal to the Lebesgue measure. In this case, it is well known that the corresponding DF process $D$ can be represented as the random probability generated on $[0, 1]$ by an increasing process of the type $\{G_t/G_1 : t \in [0, 1]\}$, where $G$ is a Gamma process on $[0, 1]$, that is, a Lévy process on $[0, 1]$ with Lévy measure (see, e.g., [15]) given by $\nu(\mathrm{d}x) = \mathbf{1}_{(x > 0)} \exp(-x)\,\mathrm{d}x/x$ (observe that the normalized process $G/G_1$ is independent of $G_1$). It follows that every $F \in L^2(D)$ is also a member of $L^2(G)$, that is, the space of square-integrable functionals of $G$. It would therefore be interesting to find some explicit relation between our orthogonal decomposition of $L^2(D)$ and the chaotic decompositions of $L^2(G)$ established in, for example, [22] or [15].

# Acknowledgements

I am indebted to Dario Spano for an introduction to Griffiths' paper [10] and for several insightful remarks. I also wish to express my gratitude to Christian Houdré, Jim Pitman and Marc Yor for inspiring discussions.

# References

[1] Abramowitz, M. and Stegun, I.A. (1964). *Handbook of Mathematical Functions with Formulas, Graphs, and Mathematical Tables.* National Bureau of Standards Applied Mathematics Series **55**. MR0167642

[2] Aldous, D.J. (1983). Exchangeability and related topics. *École d'Été de Probabilités de Saint-Flour XIII. Lecture Notes in Math.* **1117** 1–198. Berlin: Springer. MR0883646

[3] Blackwell, D. (1973). Discreteness of Ferguson selections. *Ann. Statist.* **1** 356–358. MR0348905

[4] Blackwell, D. and MacQueen, J. (1973). Ferguson distribution via Pólya urn schemes. *Ann. Statist.* **1** 353–355. MR0362614

[5] El-Dakkak, O. and Peccati, G. (2007). Hoeffding decompositions and urn sequences. Preprint.

[6] Ethier, S.N. and Kurtz, T. (1986). *Markov Processes: Characterization and Convergence.* New York: Wiley. MR0838085

[7] Ethier, S.N. and Griffiths, R.C. (1993). The transition function of a Fleming–Viot process. *Ann. Probab.* **21** 1571–1590. MR1235429

[8] Ferguson, T.S. (1973). A Bayesian analysis of some nonparametric problems. *Ann. Statist.* **1** 209–230. MR0350949

[9] Ferguson, T.S. (1974). Prior distributions on spaces of probability measures. *Ann. Statist.* **2** 615–629. MR0438568

[10] Griffiths, R.C. (1979). A transition density expansion for a multi-allele diffusion model. *Adv. in Appl. Probab.* **11** 310–325. MR0526415

[11] Koroljuk, V.S. and Borovskich, Yu.V. (1994). *Theory of U-Statistics*. Dordrecht: Kluwer Academic Publishers. MR1472486

[12] Littler, R.A. and Fackerell, E.D. (1975). Transition densities for neutral multi-allele models. *Biometrics* **31** 117–123. MR0368825

[13] Mauldin, R.D., Sudderth, W.D. and Williams, S.C. (1992). Pólya trees and random distributions. *Ann. Statist.* **20** 1203–1221. MR1186247

[14] Nualart, D. (1995). *The Malliavin Calculus and Related Topics*. New York: Springer. MR1344217

[15] Nualart, D. and Schoutens, W. (2000). Chaotic and predictable representation for Lévy processes. *Stochastic Process. Appl.* **90** 109–122. MR1787127

[16] Peccati, G. (2002). Chaos Brownien d'espace-temps, décompositions de Hoeffding et problèmes de convergence associés. Ph.D. dissertation, Univ. Paris VI.

[17] Peccati, G. (2003). Hoeffding decompositions for exchangeable sequences and chaotic representation for functionals of Dirichlet processes. *C. R. Math. Acad. Sci. Paris* **336** 845–850. MR1990026

[18] Peccati, G. (2004). Hoeffding–ANOVA decompositions for symmetric statistics of exchangeable observations. *Ann. Probab.* **32** 1796–1829. MR2073178

[19] Pitman, J. (1996). Some developments of the Blackwell–MacQueen urn scheme. In *Statistics, Probability and Game Theory: Papers in Honor of David Blackwell. IMS Lecture Notes Monogr. Ser.* **30**. Hayward, CA: IMS. MR1481784

[20] James, L.F., Lijoi, A. and Prünster, I. (2006). Conjugacy as a distinctive feature of the Dirichlet process. *Scand. J. Statist.* **33** 105–120. MR2255112

[21] Schleusner, J.W. (1974). A note on biorthogonal polynomials in two variables. *SIAM J. Math. Anal.* **5** 11–18. MR0333563

[22] Segall, A. and Kailath, T. (1976). Orthogonal functionals of independent-increment processes. *IEEE Trans. Inform. Theory* **22** 287–298. MR0413257

[23] Stroock, D.W. (1987). Homogeneous chaos revisited. *Séminaire de Probabilités XXI. Lecture Notes in Math.* **1247** 1–8. Berlin: Springer. MR0941972

[24] Tsilevitch, N., Vershik, A. and Yor, M. (2001). An infinite-dimensional analogue of the Lebesgue measure and distinguished properties of the gamma process. *J. Funct. Anal.* **185** 274–296. MR1853759